\documentclass[12pt]{article}
\usepackage{amsmath,amssymb,latexsym,graphicx}
\newtheorem{theorem}{Theorem}
\newtheorem{lemma}[theorem]{Lemma}
\newtheorem{proposition}[theorem]{Proposition}
\newtheorem{remark}[theorem]{Remark}

\def\qed{$\Box$}

\def\G{\Gamma}

\def\CC{\Bbb{C}}

\def\RR{\Bbb{R}}

\def\pf{{\bf Proof }}
\begin{document}
\title{Range descriptions for the spherical mean Radon transform
\thanks{The work of the second author was supported in part by the NSF Grants DMS-9971674,
DMS-0002195, and DMS-0604778. The work of the third author was supported in part by
the NSF Grants DMS-0200788 and DMS-0456868.} }
\author{M.~Agranovsky\\ Bar-Ilan University \and P.~Kuchment\\ Texas A\&M University \and E.~T.~Quinto\\ Tufts University}
\date{}
\maketitle

\begin{abstract}
The transform considered in the paper averages a function
supported in a ball in $\RR^n$ over all spheres centered at the
boundary of the ball. This Radon type transform arises in several
contemporary applications, e.g. in thermoacoustic tomography and
sonar and radar imaging. Range descriptions for such transforms
are important in all these areas, for instance when dealing with
incomplete data, error correction, and other issues. Four
different types of complete range descriptions are provided, some
of which also suggest inversion procedures. Necessity of three of
these (appropriately formulated) conditions holds also in general
domains, while the complete discussion of the case of general
domains would require another publication.

\end{abstract}
\tableofcontents

\section{Introduction}
The  spherical mean Radon transform, which integrates a function
over all spheres centered at points of a given set, has been
studied for quite a while in relation to PDE problems (e.g.,
\cite{Asgeirsson,CH,John}). However, there has been a recent surge
in its studies, due to demands of manifold applications. These
include, among others, the recently developed thermoacoustic and
photoacoustic tomography (e.g., \cite{AmbKuch, FPR, Kruger,
Patch}, \cite{MXW}-\cite{XWAK}), as well as radar and sonar
imaging, approximation theory, mathematical physics, and other
areas \cite{ABK, AQ, CH, Leon_Radon, John, Kuch_AMS05, LP1, LP2,
LQ, NC, Norton, Nort-Linzer}. For instance, in thermoacoustic (and
photoacoustic) tomography, the spherical mean data of an unknown
function (the radiofrequency energy absorption coefficient) is
measured by transducers, and the imaging problem is to invert that
transform (e.g., \cite{Kruger}, \cite{MXW}-\cite{YXW2}). These
applications also brought about mathematical problems that had not
been studied before. Many issues of uniqueness and stability of
reconstruction, inversion formulas, incomplete data problems,
etc., are still unresolved, in spite of a substantial body of
research available (e.g., \cite{ABK, AQ, AmbKuch, And, Den,
Leon_Radon, FHR, Faw, FPR, GGG2, Gi, Kuch_AMS05, KuchQuinto,
Kunyansky, LP1, LP2, LQ, Natt2001, Nil, NC, Norton, Nort-Linzer,
Pal, Pal05, Pal_book, Patch, Q1980, Quinto2, Qrange},
\cite{MXW}-\cite{XWAK}). In this text we address the problem which
has been recently considered for the first time
\cite{AmbKuch_range, FinRak, Patch} (see also related discussions
in \cite{Nessibi1,Nessibi2}), namely the range conditions for the
spherical mean transform. In fact, as we will mention below, in an
implicit form, a part of range conditions was already present in
\cite{LP1,LP2} and later in \cite{AQ}.

For someone coming from PDEs and mathematical physics, the range
description question might seem somewhat unusual. However, it is
well known in the areas of integral geometry and tomography that
range descriptions are of crucial theoretical and practical
importance \cite{Leon_Radon, GGG1, GGG2, GelfVil, Helg_Radon,
Natt4, Natt2001}. The ranges of Radon type transforms usually have
infinite co-dimension (e.g., in spaces of smooth functions, or in
appropriate Sobolev scales), and thus infinitely many range
conditions appear. One might wonder, what is the importance of
knowing the range conditions. The answer is that, besides their
analytic usefulness for understanding the transform, they have
been used for a variety of purposes in tomography (as well as in
radiation therapy planning \cite{CQ1,CQ2,K1}): completing
incomplete data, correcting measurement errors and hardware
imperfections, recovering unknown parameters of the medium, etc.
\cite{Hertle2, Lv, Mennes, Natt3, Natt4, Noo_half, NW, Po2, Sol3,
Sol4}. Thus, as soon as the spherical mean transform started
attracting a lot of attention, researchers started looking for
range descriptions. Some range conditions (albeit they were not
called this way) were already present in \cite{LP1,LP2} and in
\cite{AQ} (see also \cite{etti}), where the sequence of
polynomials was considered arising as moments of the spherical
Radon data. In an explicit form, these conditions were formulated
recently in \cite{Patch}. They, however, as it was discovered in
\cite{AmbKuch_range}, do not describe the range completely. A
complete set of conditions was found in the two-dimensional case
in the recent paper \cite{AmbKuch_range} and for odd dimensions
(albeit, for somewhat different transforms) in \cite{FinRak}. In
all these papers, the centers of spheres of integration (i.e., the
location of tomographic transducers) were assumed to belong to a
sphere.

In this paper, we obtain range descriptions in arbitrary dimension
for the case of centers on a sphere. Moreover, we obtain several
different range descriptions that shed new light on the meaning of
the range conditions (in particular, onto the appearance of two
seemingly different subsets of conditions).

In Section \ref{S:notions}, we introduce main notations and
preliminary facts that will be needed further on. Section
\ref{S:statement_results} contains the formulation of the main
results. Theorem \ref{T:Main_ball} provides three different types
of range descriptions. Theorem \ref{T:Main_odd} establishes that
in odd dimensions moment conditions of Theorem \ref{T:Main_ball}
can be dropped (a situation analogous to the one in
\cite{FinRak}). On the other hand, it is shown in Theorem
\ref{T:estimates} that a strengthened version of moment conditions
alone describes the range. It is shown in Theorem \ref{T:sobolev}
that the results also hold in Sobolev scale, rather than in
$C^\infty$ category. The next four sections are devoted to the
proofs of these theorems. It is noticed in Section
\ref{S:general_domains} that necessity of most of the range
conditions is in fact proven for general domains, not just for a
ball. This is described in Theorem \ref{T:Main_general}. Section
\ref{S:lemmas} contains proofs or alternative proofs of some
technical lemmas. The alternative proofs are provided, since the
authors believe that they might shed extra light onto the problem.
The final two sections contain additional remarks and discussions,
and acknowledgments.

\section{Main notions and preliminary information}\label{S:notions}

In this section we introduce main notions and notations that will
be used throughout the paper. We also remind the reader of some
known facts that we will need to use.

We will be dealing with domains in $\RR^n$. The closure of a
domain $\Omega$ is denoted by $\overline{\Omega}$ and its boundary
is $\partial \Omega$. We denote the $n-$dimensional unit ball
$B=\{x\in\RR^n|\,|x|\leq 1\}$ and the unit sphere is $S=\partial
B=\{x\in\RR^n|\,|x|=1\}$. The area of $S$ is $\omega$ and the area measure
on $S$ will be $dS$ (this notation will also be used  for the surface measure
on the boundaries of other domains).
The notation $C_0^{\infty}(B)$ stands for the class of smooth functions
with the compact support in the closed unit ball. For partial derivatives,
the notations $\frac{\partial f}{\partial t}, \partial_t f$, and $f_t$ will be used.

\subsection{Spherical means}\label{SS:means}
The main object of study in this paper is the {\em spherical mean
transform} $R$ (with centers on $S=\partial B$) that takes any
function $f\in C^\infty_0 (B)$ to
\begin{equation}\label{E:Radon}
    Rf(p,t):=\frac{1}{\omega}\int\limits_{y \in S}f(p+ty)dS(y), p\in
    S.
\end{equation}
One might wonder why we require the support of $f$ to belong to
$B$. It will be explained in Section \ref{S:remarks} that there is
not much hope for explicit range descriptions, if one allows the
support of the function to spill outside the surface $S$ of the
centers.

We will also consider the cylinder $C=B\times [0,2]$ and its
lateral boundary $S\times[0,2]$.

\begin{figure}[ht]\begin{center}
  \scalebox{0.7}{\includegraphics{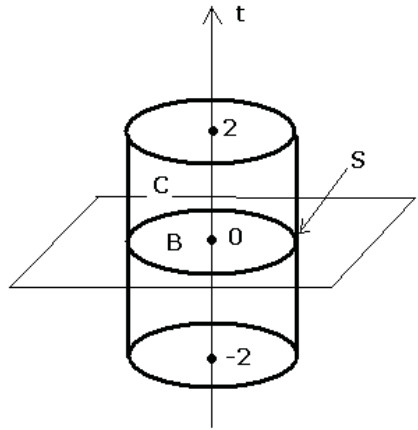}}\\
  \caption{The geometry of domains}
  \label{F:cylinder}\end{center}
\end{figure}

Notice that the sphere $S$ enters here in two different ways: as
the set of centers of the spheres of integration ($p \in S$), and
as a parametrization of the spheres of integration ($\int_{y\in
S}... dS(y)$). The reader should not confuse the two, since
sometimes in the text we will change the set of centers to a more
general surface $S$, while keeping the integration surfaces
spherical.

The results can be easily re-scaled to the case when the set of
centers is a sphere of an arbitrary radius $\rho$. We avoid doing
this here, in order to simplify the expressions.

\subsection{Darboux equation}\label{SS:darboux}

Allowing in (\ref{E:Radon}) the centers $x$ of the sphere of
integration be arbitrary, one arrives to a function
\begin{equation}\label{E:interior}
G(x,t)=\frac{1}{\omega}\int\limits_{y \in S}f(x+ty)dS(y), x\in
\RR^n.
\end{equation}
It is well known \cite{Asgeirsson,CH,John} that the function
$G(x,t)$ defined by (\ref{E:interior}) satisfies the Darboux
(Euler-Poisson-Darboux) equation
\begin{equation}\label{E:Darboux}
\frac{\partial ^2 G}{\partial t^2}+\frac{n-1}{t}\frac{\partial
G}{\partial t}=\Delta_x G,
\end{equation}
as well as the initial conditions
\begin{equation}\label{E:Darboux_init}
    G(x,0)=f(x), \frac{\partial
G}{\partial t}(x,0)=0, x\in \RR^n.
\end{equation}
Moreover, any such solution of (\ref{E:Darboux}) -
(\ref{E:Darboux_init}) in $\RR^n\times \RR^+$ is representable as
the spherical mean (\ref{E:interior}) of $f(x)$ (Asgeirsson's
Theorem, see \cite{Asgeirsson,CH,John}).

An important remark about the initial conditions
(\ref{E:Darboux_init}) is that they mean that the solution $G$ can
be extended to all values of time as an even solution on
$\RR_x^n\times\RR_t$ \cite[Ch. VI.13]{CH}.

One notices that the restriction $g$ of $G$ to $S\times \RR^+$
coincides with $Rf$:
\begin{equation}\label{E:Darboux_bound}
g=G|_{S\times \RR^+}=Rf
\end{equation}

Another observation concerning the mean $G(x,t)$ of a function
$f(x)$ supported in $B$ is that it vanishes for $x\in B, t\geq 2$.
Indeed, the value $G(x,t)$ is the average of $f$ over the sphere
centered at $x$ and of radius $t$, while such a sphere for $x\in
B, t\geq 2$ does not intersect the support of $f$, which is
contained in $B$. So, $G$ satisfies the terminal conditions
\begin{equation}\label{E:Darboux_term}
    G(x,2)=0, \frac{\partial
G}{\partial t}(x,2)=0, x\in B.
\end{equation}

We also need to mention some other known properties of the Darboux
equation, which we will need to utilize further on in the text.

The Darboux equation has a useful connection with the wave
equation (e.g., \cite[Ch. VI.13]{CH}). This relation comes from
existence of transformations that intertwine the second derivative
operator $\frac{d^2}{dt^2}$ with the Bessel operator
$\mathcal{B}_p=\frac{d^2}{dt^2}+\frac{2p+1}{t}\frac{d}{dt}$. A
general approach to constructions of such transformation operators
(not just for the Bessel case) can be found, for instance in
\cite{Levitan_shift} (see also \cite{Kipriyanov}). Among those,
the most commonly used ones are the Poisson transform (also called
Delsarte or Riemann-Liouville transform
\cite{CH,Delsarte,Levitan,Lions_dels,trimeche}) and Sonine
transform \footnote{All such transforms that we use are partial
cases of transforms of Abel type.}. The formula we will need for
the Poisson transform is
\begin{equation}
\left(\mathcal{P}_pU\right)(t)=c_p t^{-2p}\int\limits^t_0
U(y)(t^2-y^2)^{p-1/2}dy.
\end{equation}
Here $c_p$ is a non-zero constant, whose specific value is of no
relevance to us. We will use this transform for even functions
$U$, in which case it can be rewritten as
\begin{equation}\label{E:Poisson}
\left(\mathcal{P}_pU\right)(t)=c_p \int\limits^1_{-1} U(\mu
t)(1-\mu^2)^{p-1/2}d\mu.
\end{equation}
The properties of this transform are well known (e.g.,
\cite{Kipriyanov,Levitan,Levitan_shift,trimeche}). In particular,
it is known to preserve evenness of functions. We will also need
the inversion formula in the case when $2p$ is an odd integer:
\begin{equation}\label{E:Poisson_inversion}
   \left(\mathcal{P}_p^{-1}G\right)(t)=\mbox{const }t\left(\frac{d}{d(t^2)}\right)^{p+1/2}\left(t^{2p}G(t)\right).
\end{equation}

We will not use the Sonine transform here.
We, however, will be interested in the intertwining
operator sometimes called Weyl transform\footnote{This transform
is dual to a Poisson transform \cite{Kipriyanov,trimeche}.}
\cite{trimeche}:
\begin{equation}\label{E:Weyl}
U(t)=\left(\mathcal{W}_pG\right)(t):=\frac{2\Gamma
(p+1)}{\sqrt{\pi}\Gamma(p+1/2)}\int\limits^\infty_{t}
G(s)(s^2-t^2)^{p-1/2}s \,ds.
\end{equation}
which we will use for the specific value $p=(n-2)/2.$
The inverse transform $\mathcal{W}^{-1}$ for this case is given by
\begin{equation}\label{E:Weyl_inverse}
G(t)=\left(\mathcal{W}^{-1}U\right)(t)=\begin{cases}
\mbox{const} \int\limits_t^{\infty}(s^2-t^2)^{-1/2}(\frac{d}{d(s^2)})^{n/2} U(s)s ds, \mbox{ for $n$ even}\\
\mbox{const}\frac{d}{d(t^2)}^{\frac{n-1}{2}}U(t), \mbox{ for $n$
odd}.
\end{cases}
\end{equation}
One has (e.g., \cite{Kipriyanov,Levitan,trimeche}) the following
intertwining identities:
\begin{equation}\label{E:intertwine}
\begin{array}{c}
\mathcal{P}_p\frac{d^2}{dt^2}U(t)=\left(\frac{d^2}{dt^2}+\frac{2p+1}{t}\frac{d}{dt}\right)\mathcal{P}_pU(t),\\
\frac{d^2}{dt^2}\mathcal{W}_pG(t)=\mathcal{W}_p\left(\frac{d^2}{dt^2}+\frac{2p+1}{t}\frac{d}{dt}\right)G(t).
\end{array}
\end{equation}
These identities show that $G(x,t)$ (for $t\geq 0$) is a solution
of Darboux equation if and only if its Weyl transform (with
respect to $t$) $U(x,t)=\mathcal{W}_p G(x,t)$ solves the wave
equation.

It will be  important for us that the transform $\mathcal{W}$ (unlike the
Poisson transform) involves integration from $t$ to $\infty$.
Hence, when applied to a function on $\RR^+$ that vanishes for
$t>a$, it preserves this property. It is clear from the inversion
formulas that $\mathcal{W}^{-1}$ behaves the same way.

An important relation between Fourier, Fourier-Bessel, and Weyl
transforms will be presented in the next subsection.

\subsection{Bessel functions and eigenfunctions of the Bessel operator and
Laplacian}\label{SS:bessel} It will be convenient for us to use a
version of Bessel functions that is sometimes called {\em
normalized}, and sometimes {\em spherical} Bessel functions (e.g.,
\cite{Levitan,Magnus,trimech}):
\begin{equation}\label{E:modif_Bessel}
j_p (\lambda)=\frac{2^p\Gamma (p+1) J_p (\lambda)}{\lambda^p}.
\end{equation}
For $p=\frac{n-2}{2}$, according to Poisson's representations of
Bessel functions \cite{John,Magnus}, this is just the spherical
average of a plane wave in $\RR^n$.

We will also use standard expansions
\begin{equation}\label{E:Bessel_expansion}
\begin{array}{c}
J_p(\lambda)=\frac{\lambda^p}{2^p}\sum\limits_{k=0}^\infty (-1)^k \frac{\lambda ^{2k}}{2^{2k}k!\Gamma (p+k+1)},\\
j_p(\lambda)=\sum\limits_{k=0}^\infty C_k \lambda^{2k},
\end{array}
\end{equation}
where $C_k$ are non-zero constants.

The function $z(t)=j_p (\lambda t)$ satisfies the Bessel equation
\begin{equation}\label{E:Bessel}
\mathcal{B}_p z:=\frac{d^2
z}{dt^2}+\frac{2p+1}{t}\frac{dz}{dt}=-\lambda^2 z
\end{equation}
and initial conditions
$$
z (0)=1, z^\prime_\nu (0)=0.
$$

The standard Fourier-Hankel (also called Hankel, or
Fourier-Bessel) transform $g(t)\mapsto \mathcal{F}_p
(g)(s\lambda)$ and its inverse can be nicely written in terms of
$j_p$:
\begin{equation}\label{E:F-H}
\begin{array}{c}
  \mathcal{F}_p (g)(s\lambda)=\int\limits^\infty_0 g(t)j_p(\lambda t)t^{2p+1}dt \\
  g(t)=\frac{1}{2^{2p}\Gamma^2(p+1)} \int\limits^\infty_0 \mathcal{F}_p (g)(t\lambda)
  j_p(\lambda t)\lambda^{2p+1}d\lambda.\\
\end{array}
\end{equation}

We will use notation $\mathcal{F}$ for the standard
one-dimensional Fourier transform.

The following relation between the Fourier, Fourier-Hankel, and
Weyl transform (e.g., \cite[p. 124]{trimeche}) helps to understand
some parts of the further calculations:
\begin{equation}\label{E:Fourier-Weyl}
\mathcal{F}_p=c_p\mathcal{F} \mathcal{W}_p,
\end{equation}
where $c_p$ is a non-zero constant, explicit value of which is of
no relevance to this study. This relation shows that the Weyl
transform is the ratio of the Fourier and Fourier-Bessel
transforms.

We will need to use the following known Paley-Wiener theorem
\cite{Akhiezer1,Griffith,Kipriyanov,trimeche} for Fourier-Bessel
transform\footnote{Albeit the statement of the Lemma is formulated
in the works cited only for the case $H=\CC$, the proofs allow one
to consider without any change $H$-valued functions as well.}.

\begin{lemma}\label{L:PW-Bessel}
Let $p=(n-2)/2$ and $H$ be a Hilbert space. An $H$-valued function
$\Phi(\lambda)$ on $\RR$ can be represented as the transform
(\ref{E:F-H}) of an even function $g\in C^\infty_0(\RR,H)$
supported on $[-a,a]$, if and only if the following conditions are
satisfied:
\begin{enumerate}

\item $\Phi(\lambda)$ is even.

\item $\Phi(\lambda)$ extends to an entire function in $\CC$ with
Paley-Wiener estimates
\begin{equation}\label{E:PW-estimates}
\|\Phi(\lambda)\|\leq C_N (1+|\lambda|)^{-N}e^{a|\mbox{{\em Im }}
\lambda|}
\end{equation}
for any natural $N$.
\end{enumerate}
\end{lemma}
Although this result is well known, for reader's convenience, we
provide its proof in Section \ref{S:lemmas}.

We will need some special solutions of Darboux equation in the
cylinder $B\times \RR$. Let $-\lambda^2$ be in the spectrum of the
Dirichlet Laplacian in the ball $B$, and $\psi_\lambda (x)$ be the
corresponding eigenfunction, i.e.
\begin{equation}\label{E:Laplace_eigenf}
\begin{cases}
    \Delta \psi_\lambda (x)=- \lambda^2 \psi_\lambda (x) \mbox{ in }B\\
    \psi_\lambda (x)=0 \mbox{ on } S.
\end{cases}
\end{equation}
Equations (\ref{E:Bessel}) and (\ref{E:Laplace_eigenf}) imply that
the function
\begin{equation}\label{E:Darboux_eigenf}
    u_\lambda (x,t)=\psi_\lambda (x)j_{n/2-1}(\lambda t).
\end{equation}
satisfies Darboux equation (\ref{E:Darboux}).

Finally, we need descriptions of (generalized) eigenfunctions
$$
-\Delta \psi=\lambda^2 \psi.
$$
of the Laplace operator in $\RR^n$and in $B$ in terms of their
spherical harmonics expansions.

Let $Y_l^m(\theta)$, $\theta\in S^{n-1}$, $m=0, 1, ...,  1\leq l
\leq d(m)$, be an orthonormal basis of spherical harmonics, where
$m$ is the degree of the harmonic and
$$
d(m)=(2m+n-2)\frac{\G(m+n-2)}{\G(n-1)\G(m+1)}.
$$
It is known (and can be easily shown by separation of variables)
that any function
\begin{equation}\label{E:Laplace_eigenf_harm}
\begin{array}{c}
\phi_{m,l}(x)=(\lambda r)^{1-n/2}J_{n/2-1+m}(\lambda r)
Y^m_l(\theta)\\
=(\lambda r)^{m}j_{n/2-1+m}(\lambda r) Y^m_l(\theta),\mbox{ where
}x=r\theta
\end{array}
\end{equation}
is a generalized eigenfunction of the Laplace operator\footnote{As
it is customary, we use the term ``generalized eigenfunction'' for
any solution of the equation $\Delta u=-\lambda^2 u$ in $\RR^n$.}
in $\RR^n$ with the eigenvalue $-\lambda^2$. In fact, one can show
that any generalized eigenfunction of the Laplace operator in
$\RR^n$ has the following expansion into spherical harmonics:
\begin{equation}\label{E:sph_expansion}
u(r\theta)=\sum\limits_{l,m} c_{l,m} (\lambda
r)^{m}j_{n/2-1+m}(\lambda r) Y^m_l(\theta),
\end{equation}
where $r=|x|,\theta=\dfrac{x}{|x|}$. One can describe precisely
the conditions on the coefficients $c_{l,m}$, necessary and
sufficient for (\ref{E:sph_expansion}) to provide all generalized
eigenfunction, as well as generalized eigenfunctions with some prescribed growth
condition at infinity \cite{Agmon_helm1,Agmon_helm2}.

If one chooses only the values of $\lambda\neq 0$ that are zeros
of $J_{m+n/2-1}(\lambda)$, one arrives to the eigenfunctions of
the Dirichlet Laplacian in the unit ball $B$. In particular,
functions (\ref{E:Laplace_eigenf_harm}) for all $m, l$ and
$\lambda\neq 0$ such that $J_{m+n/2-1}(\lambda)=0$, form a
complete set of eigenfunctions.

With all these preparations in place, we can now set out to
formulate and prove the results of this article.

\subsection{Some preliminary results}
We start considering the spherical mean transform $R$ introduced
in (\ref{E:Radon}) with centers on the unit sphere $S$. The
following moment conditions for this case were present in
\cite{LP1,LP2}, as well as in \cite{AQ} (see also \cite{etti}),
and explicitly formulated as range conditions in \cite{Patch}:
\begin{lemma}\label{L:moment_ball}
Let $g(x,t)=Rf(x,t)$ for $f\in C^\infty_0(B)$. Then, for any
non-negative integer $k$, the function $M_k$ on $S$ defined as
\begin{equation}\label{E:Moment}
  M_k(x):=\int\limits_0^2 t^{2k+n-1} g(x,t)dt, x\in S
\end{equation}
has an extension to $\RR^n$ as a polynomial $q_k(x)$ of degree at
most $2k$.
\end{lemma}
\pf One readily observes that, for $|x|=1$,
\begin{equation}\label{E:polynomial}
M_k(x)=\int\limits_{\RR^n} |x-p|^{2k}f(p)dp.
\end{equation}
Applying (\ref{E:polynomial}) to arbitrary $x\in \RR^n$ (not
necessarily on the unit sphere), one clearly gets a polynomial
$q_k(x)$ of $x$ of degree at most (not necessarily equal to) $2k$.
\qed

\begin{remark}\label{R:moment_ball}
In fact, noticing that
$$
M_k(x):=\int\limits_{\RR^n} |x-p|^{2k}f(p)dp=\int\limits_{\RR^n}
(1-2p\cdot x+|p|^2)^{k}f(p)dp,
$$
we conclude that $M_k(x)$ has an extension to $\RR^n$ as a
polynomial of degree at most $k$.

\end{remark}

The possibility of reducing to degree $k$ comes from the fact that
we use centers of the spheres of integration that belong to the
(unit) sphere. It is clear that when the centers run over a
non-spherical surface, this reduction is not possible anymore, and
one cannot guarantee degree less than $2k$. However, in this case
an extra condition can be found that straightens up the situation.
We thus provide here an alternative reformulation of the moment
conditions, which will be handy in more general considerations.

\begin{lemma}\label{L:moment_general} Let $D\subset \RR^n$ be a
smooth bounded domain and $R$ denote the spherical mean transform
with centers on $S=\partial D$.

Let $g(x,t)=Rf(x,t)$ for $f\in C^\infty_0(D)$. Then, for any
non-negative integer $k$, the function $M_k$ on $S$ defined as in
(\ref{E:Moment}) has an extension $Q_k(x)$ to $\RR^n$ as a
polynomial of degree at most $2k$, satisfying the following
additional condition:
\begin{equation}\label{E:Moment_relation}
\Delta Q_k=c_k Q_{k-1}  \mbox{ for any }k=1,...,
\end{equation}
where $c_k=2k(2k+n-2).$
\end{lemma}
\pf Let us notice that $M_k =|x|^{2k}* f(x)|_S$. We can now define
the polynomials $Q_k(x)=|x|^{2k}* f(x)$. Then the relation follows
from the easily verifiable identity $\Delta
|x|^{2k}=c_k|x|^{2(k-1)}$. \qed

In fact, in the case when $D$ is a ball $B$, the condition
(\ref{E:Moment_relation}) is not needed.

\begin{lemma}\label{L:Moment_equiv}
In the case when $D=B$, the  condition (\ref{E:Moment_relation})
of Lemma \ref{L:moment_general} on a function $g(x,t)$ on $S\times
[0,2]$ can be dropped, and thus conditions of Lemmas
\ref{L:moment_ball} and \ref{L:moment_general} are equivalent.
\end{lemma}
\pf To prove the lemma, we need the following well known fact
(e.g., \cite{Khavinson_book}), which we prove here for the sake of
completeness.
\begin{proposition}\label{P:polyn}
The solution of the boundary value problem
$$\Delta u=v, |x|<1,$$
$$u=0, |x|=1,$$
where $v$ is a polynomial, is a polynomial of degree $\deg u=\deg
v+2.$
\end{proposition}
\pf Let us prove first that there exists a polynomial solution
$\tilde u$ of Poisson  equation
$$\Delta \tilde u=v$$ in the unit ball, such that $\deg \tilde u =\deg v+2$.
Clearly, it suffices to do this for each homogeneous term of $v$,
so we can assume the polynomial $v$ to be homogeneous.

Let us represent $v$ in the form:
$$v(x)=\sum_{\nu=0}^{[\frac{\deg v}{2}]}|x|^{2\nu}h_{\nu}(x),$$
where each $h_{\nu}$ is either zero, or a homogeneous harmonic
polynomial of degree $\deg h_{\nu}= \deg v - 2\nu$, and brackets
$[...]$ denote the integer part. This representation is well known
to be always possible (e.g., \cite{Helg_groups,Stein}). A solution
$\tilde u$ can be found in the similar form (where we denote for
brevity $k=[\frac{\deg v+2}{2}]$)
$$
\tilde u(x)=\sum_{\nu=0}^{k} |x|^{2\nu} \tilde h_{\nu}(x), \ \deg
\tilde h_{\nu} =\deg v+2-2\nu.
$$
Here again, each $\tilde h_{\nu}$ is either zero, or a homogeneous
harmonic polynomial of degree $\deg \tilde h_{\nu}=\deg v+2-2\nu$.
Then direct calculation shows
$$
\Delta \tilde u(x)=\sum_{\nu=0}^{k} \left(c_{\nu}+ 2\nu \deg
\tilde h_{\nu}\right) |x|^{2\nu-2} \tilde h_{\nu}(x),
$$
with the coefficients $c_{\nu}$ defined in Lemma 3.
Thus, the needed polynomial solution $\tilde u$ that we are
looking for can be obtained by choosing
$$
\tilde h_{\nu}=[c_{\nu} +2\nu (\deg v+2-2\nu)]^{-1} h_{\nu-1}.
$$
To finish the proof of Proposition, introduce $U=\tilde u -u$.
Then one obtains
$$\Delta U=0, |x|<1,$$
$$U=\tilde u, |x|=1,$$
for the newly defined function $U$. The boundary value $\tilde u$ is the
polynomial of degree $\deg \tilde u=\deg v+2$. Its harmonic
extension $U$ from the unit sphere is obtained from the above
decomposition of $\tilde u$ by replacing $|x|$ by $1$:
$$
U(x)=\sum\limits_{\nu =0}^k \tilde h_\nu (x).
$$
Since $\deg \tilde h_\nu \leq \deg v +2$, we have $\deg U \leq
\deg v+2$.

Thus, the solution $u=\tilde u-U$ is a polynomial of degree at
most $\deg v+2$, which proves the Proposition. \qed

\begin{remark}\label{R:polyn}
Polynomial solvability of the Poisson problem with polynomial data
is rather unique and essentially holds only for balls (e.g.,
\cite{Siegel,shapiro,Khavinson}). Hence, in a general domain, the
conditions of Lemma \ref{L:moment_general} are stronger than the
ones of Lemma \ref{L:moment_ball}.
\end{remark}

{\bf Proof of Lemma \ref{L:Moment_equiv}.} We want to prove that
among all polynomial extensions $q_k$, $\deg q_k \le 2k,$ of the
functions $M_k$ defined in (17), there is a sequence of extensions
$Q_k$ satisfying the additional recurrence relation
(\ref{E:Moment_relation}). Let $q_k$ be some extensions. Any other
sequence $Q_k$ of extensions can be represented as
$$Q_k=q_k+u_k,$$
where $u_k=0$ on the unit sphere. The additional requirement
(\ref{E:Moment_relation}) yields the relation
$$
\Delta u_k=c_kQ_{k-1}-\Delta q_k.
$$
The existence of polynomial solutions $u_k, \deg u_k \le 2k$
follows now by inductive application of Proposition \ref{P:polyn}.
Then the modified sequence of polynomials $Q_k=q_k+u_k$ satisfies
all the requirements of the lemma.
 \qed

\begin{remark}
The sequence $Q_k$ of polynomial extensions of the functions $M_k$
satisfying the chain relation (19) is unique. Indeed, if
$Q_k^{\prime}$ is another sequence of such extensions, then the
polynomials $R_k=Q_k-Q_k^{\prime}$  vanish on the unit sphere and
still possess (\ref{E:Moment_relation}). If not all $R_k$ are
identically zero, then, due to (\ref{E:Moment_relation}), the
first nonzero polynomial $R_{k_0}$ is harmonic and vanishes on the
unit sphere. Thus, it must be zero, due to the maximum principle.
This contradiction shows that $R_k=0$ and hence $Q_k=Q_k^{\prime}$
for all $k$.
\end{remark}

In the future, we  will need the
moment condition on the unit sphere in a different form
\cite{AmbKuch_range}:
\begin{lemma}\label{L:moments_interpretation}
The function $g(x,t),  x \in S,  t \in [0,2]$ satisfies
(\ref{E:Moment}) if and only if, for any spherical harmonic
$Y^m(\theta)$ of degree $m$, the function
$$
\hat g_m(\lambda)= \int\limits_0^2 \int\limits_{\theta \in S}
g(\theta,t) j_{n/2-1+m}(\lambda t)Y^{m}(\theta)dS(\theta)
t^{n-1}dt
$$
has at $\lambda=0$ a zero of order at least $m$.
\end{lemma}
\pf The moment conditions require that
$$
\int\limits_0^\infty t^{2k+n-1}g(\theta,t)dt
$$
is extendable to a polynomial of degree at most $k$. Let us expand $g(\theta,t)$ into an
orthonormal basis $Y^m_l$ of spherical harmonics on $S$ (where $m$ is the degree of the harmonic):
$$
g(\theta,t)=\sum\limits_{l,m} g_{l,m}(t)Y^m_l(\theta).
$$
Due to smoothness and compactness of support of $g$, it is legitimate to integrate term-wise in computing the momenta to obtain
\begin{equation}\label{E:momenta_harmonics}
    \int\limits_0^\infty t^{2k+n-1}g(\theta,t)dt=\sum\limits_{l,m} Y^m_l(\theta)\int\limits_0^\infty t^{2k+n-1}g_{l,m}(t)dt.
\end{equation}
A spherical harmonic of degree $m$ can be extended to a polynomial
of degree $d$ if and only  if $d\geq m$ and $d-m$ is even. Thus,
the moment conditions require that the coefficients $\int_0^\infty
t^{2k+n-1}g_{l,m}(t)dt$ in (\ref{E:momenta_harmonics}) with $m>2k$
must vanish.

Let us turn to the function $\hat g_m(\lambda)$.
Using the expansion
(\ref{E:Bessel_expansion}) for Bessel functions, one arrives to the formula
\begin{equation}\label{E:g_hat_expansion}
\hat g_m (\lambda)=\sum\limits_{k=0}^\infty C_k\lambda^{2k}
\int\limits_0^\infty t^{2k+n-1}\int\limits_{\theta \in
S}Y^m(\theta)g(\theta,t)dS(\theta) dt .
\end{equation}
Now one sees that the moment conditions are equivalent to the
requirement that all terms in this series with $2k<m$ vanish.
Therefore, the series begins with the power at least
$\lambda^{m}.$ \qed

\section{Statements of the main results}\label{S:statement_results}

As it follows from \cite{AmbKuch_range,FinRak}, the moment range
conditions of the preceding lemmas are insufficient. Necessary and
sufficient conditions of different kinds were provided in
\cite{AmbKuch_range} in dimension two and in \cite{FinRak} in odd
dimensions (albeit for somewhat different transforms). We
formulate below our main result that resolves the problem of range
description in any dimension, as well as provides several
alternative ways to describe the range.

\begin{theorem}\label{T:Main_ball}
    The following four statements are equivalent:
    \begin{enumerate}
    \item The function $g\in C^\infty_0 (S\times [0,2])$ is
    representable as $Rf$ for some $f\in C^\infty_0(B)$.
    \item
    \begin{enumerate}
    \item The moment conditions of Lemma \ref{L:moment_ball} (or Lemma \ref{L:moment_general}) are satisfied.
    \item The solution $G(x,t)$ of the interior problem
    (\ref{E:Darboux}), (\ref{E:Darboux_bound}), (\ref{E:Darboux_term}) in $C$ (which always exists for $t>0$) satisfies the condition
    $$
    \lim\limits_{t \to 0}\int\limits_B \frac{\partial G}{\partial t}(x,t)\phi(x)dx=0
    $$
    for any eigenfunction $\phi(x)$ of the Dirichlet Laplacian in
    $B$.
    \end{enumerate}
    \item
    \begin{enumerate}
    \item The moment conditions of Lemma \ref{L:moment_ball} (or Lemma \ref{L:moment_general}) are satisfied.
    \item Let $-\lambda^2$ be an eigenvalue of the Dirichlet Laplacian in $B$
    and $u_\lambda$ be the corresponding eigenfunction solution (\ref{E:Darboux_eigenf}).
    Then the following orthogonality condition is satisfied:
    \begin{equation}\label{E:orthog}
\int\limits_{S\times [0,2]} g(x,t)\partial_\nu u_\lambda
(x,t)t^{n-1}dxdt=0.
    \end{equation}
    Here $\partial_\nu$ is the exterior normal derivative at the
    boundary of $C$.
    \end{enumerate}
    \item
    \begin{enumerate}
    \item The moment conditions of Lemma \ref{L:moment_ball} (or Lemma \ref{L:moment_general}) are satisfied.
    \item Let $\widehat{g}(x,\lambda)=\int g(x,t)j_{n/2-1}(\lambda t)t^{n-1}dt$. Then, for any integer $m$, the
    $m^{th}$ order spherical harmonic term $\widehat{g}_m(x,\lambda)$ of $\widehat{g}(x,\lambda)$ vanishes at
    non-zero zeros of the Bessel function $J_{m+n/2-1}(\lambda)$.
    \end{enumerate}
    \end{enumerate}
\end{theorem}

In fact, if dimension $n$ is odd, one does not need to require the
moment conditions of Lemma \ref{L:moment_ball} or Lemma
\ref{L:moment_general}, since in this case the other range
conditions of Theorem \ref{T:Main_ball} alone are sufficient. This
was first noticed in a different setting in \cite{FinRak}.

\begin{theorem}\label{T:Main_odd} Let $n>1$ be an odd integer. Then the following four statements are equivalent:
    \begin{enumerate}
    \item The function $g\in C^\infty_0 (S\times [0,2])$ is
    representable as $Rf$ for some $f\in C^\infty_0(B)$.
    \item
    The solution $G(x,t)$ of the interior problem
    (\ref{E:Darboux}), (\ref{E:Darboux_bound}), (\ref{E:Darboux_term}) in $C$ (which always exists for $t>0$) satisfies the condition
    $$
    \lim\limits_{t \to 0}\int\limits_B \frac{\partial G}{\partial t}(x,t)\phi(x)dx=0
    $$
    for any eigenfunction $\phi(x)$ of the Dirichlet Laplacian in
    $B$.
    \item
    Let $-\lambda^2$ be an eigenvalue of the Dirichlet Laplacian in $B$
    and $u_\lambda$ be the corresponding eigenfunction solution (\ref{E:Darboux_eigenf}).
    Then the following orthogonality condition is satisfied:
    \begin{equation}
\int\limits_{S\times [0,2]} g(x,t)\partial_\nu u_\lambda
(x,t)t^{n-1}dxdt=0.
    \end{equation}
    Here $\partial_\nu$ is the exterior normal derivative at the
    boundary of $C$.
    \item
    Let $\widehat{g}(x,\lambda)=\int g(x,t)j_{n/2-1}(\lambda t)t^{n-1}dt$. Then, for any integer $m$, the
    $m$th order spherical harmonic term $\widehat{g}_m(x,\lambda)$ of $\widehat{g}(x,\lambda)$ vanishes at
    non-zero zeros of the Bessel function $J_{m+n/2-1}(\lambda)$.
    \end{enumerate}
\end{theorem}

As we proved in Lemma \ref{L:Moment_equiv}, in the case of a ball
the polynomial extendibility of the $t$-moments $M_k$ of the data
$g(x,t)$ (the moment condition) is equivalent to the stronger
moment condition (Lemma \ref{L:moment_general}) that requires
existence of polynomial extensions $Q_k$ linked by the additional
relations (\ref{E:Moment_relation}).

In fact, if $g$ is in the range $g=Rf, f \in C_0(B)$, then the
polynomials $Q_k$ obey not only algebraic relations
(\ref{E:Moment_relation}), but also the following growth (in $k$)
estimates in $B$:
\begin{equation}\label{E:Q_estimate}
|Q_k(x)|=|\int\limits_B|x-y|^{2k}f(y)dy| \le 2^{2k} \max\limits_{y
\in B}|f(y)|, x \in B.
\end{equation}

It turns out that the moment conditions of Lemma
\ref{L:moment_general} with estimates of the above type are not
only necessary, but also sufficient for $g$ being in the range of
the transform $R$. The following theorem can be regarded as an
alternative form of Theorem \ref{T:Main_ball}:

\begin{theorem}\label{T:estimates}
Let $g \in C_0^{\infty}(S \times [0,2])$. The following condition
is necessary and sufficient for the function $g$ being
representable as  $g=Rf$ for some $f \in C_0^{\infty}(B)$:

The moments
$$
M_k(x)=\int\limits_0^{\infty}g(x,t)t^{2k+n-1}dt
$$
extend from $x \in S$ to $x \in \mathbb R^n$ as polynomials
$Q_k(x)$ satisfying the recurrent condition
(\ref{E:Moment_relation}) and the growth estimates
\begin{equation}\label{E:ball_estimate}
|Q_k(x)|\le M^k, x \in B,
\end{equation}
for some $M>0$.
\end{theorem}

The condition of infinite smoothness of functions under
consideration is not truly necessary. One can prove similar range
descriptions in appropriate Sobolev spaces, if the functions are
supported strictly inside the ball $B$.
\begin{theorem}\label{T:sobolev}
Let $s\geq 0$. The following four statements are equivalent:
    \begin{enumerate}
    \item The function $g\in H_{s+(n-1)/2}^{\mbox{comp}} (S\times (0,2))$ is
    representable as $Rf$ for some $f\in H_s^{\mbox{comp}}(B)$.
    \item
    \begin{enumerate}
    \item The moment conditions of Lemma \ref{L:moment_ball} (or Lemma \ref{L:moment_general}) are satisfied.
    \item The solution $G(x,t)$ of the interior problem
    (\ref{E:Darboux}), (\ref{E:Darboux_bound}), (\ref{E:Darboux_term}) in $C$ (which always exists for $t>0$) satisfies the condition
    $$
    \lim\limits_{t \to 0}\int\limits_B \frac{\partial G}{\partial t}(x,t)\phi(x)dx=0
    $$
    for any eigenfunction $\phi(x)$ of the Dirichlet Laplacian in
    $B$.
    \end{enumerate}
    \item
    \begin{enumerate}
    \item The moment conditions of Lemma \ref{L:moment_ball} (or Lemma \ref{L:moment_general}) are satisfied.
    \item Let $-\lambda^2$ be an eigenvalue of the Dirichlet Laplacian in $B$
    and $u_\lambda$ be the corresponding eigenfunction solution (\ref{E:Darboux_eigenf}).
    Then the following orthogonality condition is satisfied:
    \begin{equation}\label{E:orthog}
\int\limits_{S\times [0,2]} g(x,t)\partial_\nu u_\lambda
(x,t)t^{n-1}dxdt=0.
    \end{equation}
    Here $\partial_\nu$ is the exterior normal derivative at the
    boundary of $C$.
    \end{enumerate}
    \item
    \begin{enumerate}
    \item The moment conditions of Lemma \ref{L:moment_ball} (or Lemma \ref{L:moment_general}) are satisfied.
    \item Let $\widehat{g}(x,\lambda)=\int g(x,t)j_{n/2-1}(\lambda t)t^{n-1}dt$. Then, for any integer $m$, the
    $m^{th}$ order spherical harmonic term $\widehat{g}_m(x,\lambda)$ of $\widehat{g}(x,\lambda)$ vanishes at
    non-zero zeros of the Bessel function $J_{m+n/2-1}(\lambda)$.
    \end{enumerate}
    \end{enumerate}

    When $n$ is odd, the moment conditions can be dropped.
\end{theorem}
Here we used the notation $H_s^{\mbox{comp}}(B)$ for the space of
$H_s$-functions in the ball $B$ with compact support in the open
ball. Analogously, $H_{s}^{\mbox{comp}} (S\times (0,2))$ consists
of $H_s$-functions on $S\times (0,2)$ with support in $S\times
(\varepsilon,2-\varepsilon)$ for some positive $\varepsilon$.

\section{Proof of Theorem \ref{T:Main_ball}}\label{S:proof_ball}
\subsection{Implication $1\Rightarrow 2$}
Assume that $1)$ is satisfied, i.e. $g(p,r)=Rf(p,r)$ for a smooth
function $f$ supported in $B$.

The implication $1\Rightarrow 2(a)$ is the statement of Lemma
\ref{L:moment_ball}.

The implication $1\Rightarrow 2(b)$ is one of the statements of
Asgeirsson's theorem \cite{Asgeirsson,CH,John,Helg_Radon}, which
has already been quoted before.

\subsection{Equivalence $2 \Leftrightarrow 3$}
Since conditions $2(a)$ and $3(a)$ are the same, we only need to
establish the equivalence $2(b) \Leftrightarrow 3 (b)$. This is
done in the lemma below. Notice that this lemma applies to any
bounded domain, not just to a ball.

\begin{lemma}\label{L:equiv}
Let $D$ be a bounded domain in $\mathbb R^n$ with smooth boundary,
$\Delta_D$ - the Laplacian in $D$ with the Dirichlet boundary
conditions, and $T>0$ be such that all spheres of radius $T$
centered in $\overline{D}$ do not intersect $\overline{D}$. Let
also $g\in C^\infty_0 (\partial D \times [0,T])$.

\begin{enumerate}
\item For any eigenfunction $\phi(x)=\phi_k(x)$ of $\Delta_D$ with
the eigenvalue $-\lambda^2=-\lambda_k^2 \in \sigma(\Delta_D)$, the
following two statements are equivalent:
\begin{enumerate}
\item
$$
\int\limits_0^{\infty} \int\limits_{\partial D}
g(x,t)\partial_\nu u_{\lambda}(x,t)t^{n-1}ds(x)dt=0,
$$
where
$$
u_{\lambda}(x,t)=\phi(x)j_{n/2-1}(\lambda t).
$$

\item The solution to the backward initial value boundary value problem
\begin{equation}\label{E:reversed}
\begin{cases}
\dfrac{\partial^2 G}{\partial t^2}+\dfrac{n-1}{t}\dfrac{\partial
G}{\partial t}=\Delta_x G, \ (x,t) \in D \times
(0,T]\\
G(x,T)=0, \partial_t G(x,T)=0; \\ G(x,t)=g(x,t) \ (x \in \partial
D)
\end{cases}
\end{equation}
(which always exists) satisfies the condition $\int\limits_D
\partial_t G(x,t)\phi (x) dx \to 0$ as $t \to 0+$.
\end{enumerate}

\item If the equivalent conditions (a) and (b) hold for all
Dirichlet eigenfunctions $\phi=\phi_k$, then there exists a smooth
function $f(x)$ in $D$, such that $\lim\limits_{t \to 0+}
G(x,t)=f(x)$, and $\lim\limits_{t \to 0+} G_t(x,t)=0$, where the
limits are understood as convergence in $C^\infty(D)$ (or,
equivalently, in the Sobolev space $H^s(D)$ for arbitrary $s$).

The converse statement also holds, i.e. this behavior of $G$ at
$t\to 0+$ implies (a) and (b) for any eigenfunction $\phi$.
\end{enumerate}
\end{lemma}

\pf of the lemma.

{\bf 1. Proof of equivalence of conditions 1(a) and 1(b).}

First of all, we need to be sure that a solution $G(x,t)$ of the
problem (\ref{E:reversed}) exists and is unique and regular for
$t>0$. This is immediate, due to the hyperbolic nature of this
problem (at least, until one approaches the singularity at $t=0$).
One can also show this as follows. Applying the Weyl transform
with respect to time to the functions $G(x,t)$ and $g(x,t)$ in
(\ref{E:reversed}), one arrives (as we have discussed already) to
a similar problem for the wave equation, where the corresponding
theorems are available in PDE textbooks (e.g., \cite[Section 7.2,
Theorem 6]{Evans}). Then, applying the inverse Weyl transform, one
obtains the needed solution of (\ref{E:reversed})\footnote{See,
e.g., \cite{Kipriyanov} for usage of such transformation
techniques for various PDE problems.}. In fact, a more elaborate
consideration of this kind can be found further on in this proof.

We will now prove the implication (a) $\rightarrow$ (b). We choose
a small $\epsilon
>0$ and start with a straightforward equality
$$
\int\limits^T_\epsilon \int\limits_{\partial D} g\partial_\nu
u_\lambda t^{n-1}dS dt =\int\limits^T_{\epsilon} j_{n/2-1}(\lambda
t)t^{n-1}dt\int\limits_{\partial D} g\partial_\nu \phi ds.
$$
Using Stokes' formula, one rewrites the inner integral as
\begin{equation}\label{E:stokes}
\int\limits_{\partial D} g\partial_\nu \phi
dS=\int\limits_{\partial D}\partial_\nu g \phi dS+\int\limits_D
(G\Delta \phi - \Delta G \phi)dx.
\end{equation}
The first integral on the right is zero, since $\phi$ vanishes on
the boundary. Then in the second integral, we use the
eigenfunction property for $\phi$ to get
$$
\int\limits_{\partial D} g\partial_\nu \phi dS=-\int\limits_D
(\Delta G +\lambda^2 G)\phi dx.
$$
Since $G$ satisfies Darboux equation, we can replace $\Delta G$ by
the Bessel operator $\mathcal{B}:=
\mathcal{B}_{(n-2)/2}=\frac{\partial^2}{\partial
t^2}+\frac{n-1}{t}\frac{\partial}{\partial t}$. This leads to the
following form of the last expression:
$$
-\int\limits_D (G_{tt} +(n-1)t^{-1}G_t+\lambda^2 G)\phi dx.
$$
Substituting this into the right hand side of (\ref{E:stokes}) and
changing the order of integration, one arrives to
\begin{equation}\label{E:stokes2}
\int\limits^T_\epsilon \int\limits_{\partial D} g\partial_\nu
u_\lambda dS dt=-\int\limits_D \phi(x)dx\int\limits^T_\epsilon
(G_{tt} +\frac{n-1}{t}G_t+\lambda^2 G)j_{n/2-1}(\lambda
t)t^{n-1}dt.
\end{equation}

Let us introduce a temporary notation
$$
h(t):=\int\limits_D G(x,t)\phi(x)dx \mbox{ for }t>0.
$$
Integrating by parts with respect to $t$ in the inner integral in
(\ref{E:stokes2}), we can rewrite the resulting expression for
$\int\limits^T_\epsilon \int\limits_{\partial D} g\partial_\nu
u_\lambda ds dt$ as follows:
\begin{equation}\label{E:stokes3}
\int\limits^T_\epsilon \int\limits_{\partial D} g\partial_\nu
u_\lambda ds dt=\epsilon^{n-1}\left(
h_t(\epsilon)j_{n/2-1}(\epsilon t) -h(\epsilon)(j_{n/2-1}(\lambda
t))_t|_{t=\epsilon}\right).
\end{equation}

We now need to investigate possible behavior of $h(t)$ and its
derivative when $t \to 0$. In order to do so, let us derive from
the Darboux equation for $G$ a differential equation for $h(t)$.
Applying the Bessel operator $\partial ^2 /\partial t^2
+(n-1)t^{-1}\partial /\partial t$ to the identity defining the
function $h(t)$, one obtains
\begin{equation}\label{E:h_equation}
\frac{\partial ^2 h}{\partial t^2} +\frac{n-1}{t}\frac{\partial
h}{\partial t}=\int\limits_D \left(\frac{\partial ^2 G}{\partial
t^2} +\frac{n-1}{t}\frac{\partial G}{\partial t} \right)\phi dx
=-\lambda^2 h -\int\limits_{\partial D} g\partial_\nu \phi dS.
\end{equation}
We used here the Darboux equation for $G$, integration by parts,
the fact that $\phi$ is an eigenfunction, and finally the
vanishing of $\phi$ at $\partial D$.

Let us introduce a shorthand notation for the last integral in
(\ref{E:h_equation}):
$$
w(t)=-\int\limits_{\partial D} g(x,t)\partial_\nu \phi (x) dS(x).
$$
Thus, we get the final non-homogeneous Bessel ODE for $h(t)$:
\begin{equation}\label{E:h-bessel}
\frac{\partial ^2 h}{\partial t^2} +\frac{n-1}{t}\frac{\partial
h}{\partial t}+\lambda^2 h=w(t).
\end{equation}
Due to the condition that $g$ belongs to $C^\infty_0 (\partial
D\times [0,T])$, we conclude that $w(t)$ is smooth and vanishes to
the infinite order at $t=0$. It is a matter of simple
consideration to show existence of a particular solution of
(\ref{E:h-bessel}) that vanishes to the infinite order at the
origin. Thus, the type of the behavior at the origin is dictated
by the solutions of the homogeneous equation. This behavior is
well known (e.g., \cite{Levitan}). It depends on whether Bessel
functions of the first or the second kind are involved.

Bessel functions of the first kind are smooth and have zero
derivative at the origin. The ones of the second type, have
singularity at zero. If there are no Bessel functions of the
second kind involved, then the solution of the homogeneous
equation is continuous at the origin and has zero derivative
there. On the other hand, if there is a Bessel function of the
second kind as a part of $h(t)$, then when $t \to 0+,$ $h$ behaves
as follows (e.g., \cite{Levitan}): when $n=2$, then $h(t)=\log
t(C+o(t))$ and $\ h^{\prime}(t)=t^{-1}(C+o(t))$ with non-zero
constants $C$. In the case when $n>2$, the corresponding behavior
is $h(t)=t^{2-n}(C+o(t))$ and $\ h^{\prime}(t)=t^{1-n}(C+o(t))$.
We will now show that this type of behavior is impossible, due to
(a). Indeed, (a) says that $\int^T_\epsilon \int_{\partial D}
g\partial_\nu u_\lambda ds dt \to 0$ when $\epsilon \to 0.$ Then,
due to (\ref{E:stokes3}),
$$
\epsilon^{n-1}\left( h_t(\epsilon)j_{n/2-1}(\epsilon t)
-h(\epsilon)(j_{n/2-1}(\lambda t))_t|_{t=\epsilon}\right) \to 0,
t\to 0.
$$
On the other hand, if Bessel functions of the second kind were
involved, then this expression would be $C+o(1)$ with a non-zero
constant $C$, which is a contradiction. Thus, we conclude that
$h(t)$ is continuous at $t=0$, and also that $h^\prime (0)=0$. The
latter statement is exactly the claim of $(b)$.

\begin{remark}
In fact, we have proven more than we claimed in (b). Indeed, we
showed not only that $h_t(t)=\int_D \partial_t G(x,t)\phi (x) dx
\to 0$ as $t \to 0+$, but also that $h(t)=\int_D G(x,t)\phi (x)
dx$ is continuous at $t=0$.
\end{remark}

The converse implication $(b)\rightarrow (a)$ is even simpler.
Condition (b) means that $h_t(\epsilon) \to 0, \epsilon \to 0+.$
Therefore, $h(\epsilon)$ has no singularity at $\epsilon=0$ and is
continuous there. Then the right hand side in (\ref{E:stokes3})
tends to zero as $\epsilon \to 0$ and therefore the left hand side
does as well. This means that (a) holds.

{\bf 2. Proof of statement 2: Regularity of  $G$ at $t=0$.} In
this part of the proof, we will use the transformation technique
already briefly mentioned above, which allows one to toggle
between the solutions of the wave equation and Darboux equation.
Ideologically, what we are about to do, is using the Weyl
transform. This can be done, and has been done by the authors.
However, it seemed to the authors, that using only Fourier and
Fourier -Bessel transforms makes the proof less technical and more
transparent. An alternative version of the proof, which uses Weyl
transform explicitly is provided in Section \ref{S:lemmas}.

First of all, it is well known (e.g., \cite[Ch.6.13]{CH}) that
existence of the limit when $t\to 0$ of $G(x,t)$ and the equality
$\lim_{t\to 0}G_t=0$ (even in weaker topologies than $C^\infty$)
mean that $G(x,t)$ can be extended to an even with respect to $t$
solution of the Darboux equation. Due to the zero conditions at
$t=T$, this even solution will be supported in $D\times [-T,T]$.

Thus, our task, instead of studying the limits of $G$ and $G_t$
when $t\to 0$ (which we would need to do if using Weyl transform),
will be to investigate existence and regularity (as a function of
$t$ with values in $H^s(D)$) of an even with respect to time
solution $G$.

Suppose we do have such a solution $G(x,t)$ supported in $D\times
[-T,T]$. Let us then take its Fourier-Bessel transform
$\mathcal{F}_p$ (with $p=(n-2)/2$) with respect to time
(\ref{E:F-H}), to get a function $\widehat G (x, \lambda)$.
According to the Lemma \ref{L:PW-Bessel}, this function, as an
$H^s(D)$-valued function of $\lambda$, would be even with respect
to $\lambda$ and would satisfy the Paley-Wiener estimate
(\ref{E:PW-estimates}) with $a=T$. Besides, the Darboux equation
and the boundary conditions would also imply that the following
equation and boundary conditions are satisfied:
\begin{equation}\label{E:Darboux_transformed}
\begin{cases}
(\lambda^2+\Delta_x)\widehat G(x,\lambda)=0 \mbox{ in } D\\
\widehat G(x,\lambda)=\widehat g(x,\lambda) \mbox{ for } x\in S.
\end{cases}
\end{equation}
Notice, that $\widehat g(x,\lambda)$ is even with respect to
$\lambda$ and of the appropriate Pale-Wiener class as a
$H^s(S)$-valued function for any $s$, due to our conditions on
smoothness and support of $g(x,t)$.

So, our problem is now equivalently reformulated as showing
existence of an even and entire with respect to $\lambda$ solution
of (\ref{E:Darboux_transformed}) of the appropriate Paley-Wiener
class.

It is clear that (\ref{E:Darboux_transformed}) might not have any
solution at all when $-\lambda^2$ belongs to the spectrum of the
Dirichlet Laplacian in $D$. However, the necessary and sufficient
condition for solvability of (\ref{E:Darboux_transformed}) for
such values of $\lambda$ are well known and easy to derive (they
represent the Fredholm alternative):
\begin{equation}\label{E:fredholm}
   \int\limits_S \widehat g(x,\lambda) \partial_\nu\phi (x)dS=0
\end{equation}
for any eigenfunction $\phi$ of $\Delta_D$ corresponding to the
eigenvalue $-\lambda^2$. These conditions clearly are equivalent
to $3(b)$ and thus satisfied in our case. Hence, one hopes to
solve (\ref{E:Darboux_transformed}) for all $\lambda$ and to
eventually get the needed solution. This is exactly what we will
endeavor now.

First of all, it will be convenient for us to apply the standard
trick of moving the inhomogeneity in (\ref{E:Darboux_transformed})
from the boundary condition into the equation. Let us denote by
$E$ any ``nice'' extension operator of functions from $S$ to $D$,
for instance any one that would map Sobolev spaces $H^{s}(S)$ to
$H^{s+1/2}(D)$. Existence of such operators is well known (see,
e.g., \cite{Lions}). The Poisson operator of harmonic extension is
one of them. Let us denote $U(x,\lambda)=\widehat G
(x,\lambda)-E\hat g(x,\lambda)$. Then this function solves the
problem
\begin{equation}\label{E:Darboux_inhom}
\begin{cases}
(\lambda^2+\Delta_x)U(x,\lambda)=\hat f(x,\lambda) \mbox{ in } B\\
U(x,\lambda)=0 \mbox{ for } x\in S.
\end{cases}
\end{equation}
Here $\hat f(x,\lambda)=-(\lambda^2+\Delta_x)E\hat g(x,\lambda))$
is of the same Paley-Wiener class with respect to $\lambda$, as
$\hat g(x,\lambda)$.

Let us apply to (\ref{E:Darboux_inhom}) the inverse Fourier
(rather than Fourier-Bessel) transform with respect to $\lambda$
(this amounts to applying the Weyl transform to the original
Darboux equation). Then we arrive to the following evolution
problem:
\begin{equation}\label{E:wave}
    \begin{cases}
    U_{tt}(x,t)=\Delta U(x,t) +f(x,t), x\in B, t\in \RR\\
    U(x,t)|_{x\in S}=0.
    \end{cases}
\end{equation}
Here $U$ and $f$ are inverse Fourier transforms from $\lambda$ to
$t$ of $\widehat U$ and $\hat f$. Function $f$ is even with
respect to $t$, infinitely smooth as $H^s(D)$-valued function of
$t$ for any $s$, and is supported (due to the Paley-Wiener
theorem) in $D\times [-T,T]$. Our goal now boils down to proving
existence of an even with respect to time solution $U(x,t)$ that
is smooth as $H^s(D)$-valued function of $t$ and is supported in
$D\times [-T,T]$. If this is done, then taking Fourier transform
with respect to time first and the inverse Fourier-Bessel
transform next, we will arrive to the solution $G(x,t)$ we need,
which will finish the proof of the Lemma.

Let us consider the following problem:
\begin{equation}\label{E:wave_init}
    \begin{cases}
    U_{tt}(x,t)=\Delta U(x,t) +f(x,t), x\in B, t>-T\\
    U(x,t)|_{x\in S}=0\\
    U(x,-T)=U_t(x,-T)=0.
    \end{cases}
\end{equation}
According to the standard existence and uniqueness theorems for
the wave equation (e.g., \cite[Section 7.2, Theorem 6]{Evans}),
there exists unique (and smooth as $H^s(D)$-valued function of
$t$) solution of this problem. Due to the type of the initial and
boundary conditions we imposed, one can extend the solution to all
times by assuming that it is zero for $t<-T$.

It only remains to prove that $U(x,t)$ vanishes for $t>T$ and is
even with respect to time. To do so, let us consider a complete
orthonormal set $\{\phi_k (x)\}$ of eigenfunctions of the
Dirichlet Laplacian in $D$ and denote by $-\lambda_k^2$ the
corresponding eigenvalues. Let us also expand the functions $U$
and $f$ into this basis:
\begin{equation}\label{E:eigen_expansion}
  \begin{cases}
  U(x,t)=\sum u_k(t)\phi_k(x)\\
f(x,t)=\sum f_k(t)\phi_k(x).
  \end{cases}
\end{equation}
It will be sufficient for our purpose to show that all functions
$u_k(t)$ vanish for $t>T$.

Let us notice that the following initial value problem is
satisfied by $u_k(t)$:
\begin{equation}\label{E:ODE}
\begin{cases}
 u_k^{\prime \prime}+\lambda_k^2u_k(t)=f_k(t)\\
 u_k(t)=u_k^\prime(t)=0 \mbox{ for }t\leq -2.
 \end{cases}
\end{equation}
Taking Fourier transform (in distribution sense) in (\ref{E:ODE}),
we get
\begin{equation}\label{E:ODE_Fourier}
 (\lambda_k^2-\lambda^2)\hat u_k(\lambda)=\hat f_k(\lambda).
\end{equation}
Here $\hat f_k(\lambda)$ is even and from the Paley-Wiener class
corresponding to smooth functions with support in $[-T,T]$.
Consider the function
\begin{equation}\label{E:aux_ODE}
 \hat v_k(\lambda):=\frac{\hat
 f_k(\lambda)}{(\lambda_k^2-\lambda^2)}.
\end{equation}
As it was mentioned above in this proof, conditions $3(b)$
guarantee that $\hat
 f_k(\lambda)$ vanishes at the points $\pm \lambda_k$. Thus, $\hat v_k$ is entire, even, and by
 a simple estimate, belongs to the same Paley-Wiener
 class as $\hat f_k$. Thus, it is Fourier transform of a smooth
 even function $v_k(t)$ supported in $[-T,T]$. Consider the
 difference $w_k(t)=u_k(t)-v_k(t)$. It satisfies then the homogeneous
 equation $w_k^{\prime \prime}+\lambda_k^2w_k=0$ and zero initial
 conditions at $t=-T$. Thus, it is identically zero. Hence,
 $u_k=v_k$ is even and supported in $[-T,T]$ for any $k$, and thus
 $u(x,t)$ is also even and supported in $D\times [-T,T]$.
 This finishes the proof of the existence of a solution $G_+(x,t)$ of the
Darboux equation inside the cylinder $C$ that agrees with the
spherical mean data on $S \times \RR$ and which is even with
respect to time, smooth as an $H^s(D)$-valued function of $t$, and
supported in $t \in [-T,T]$.

What now remains to prove in the lemma, is the converse statement
in its part 2. This is, however, trivial. Indeed, the strong
convergence of $G$ at $t\to 0$ we have derived clearly implies the
statement 1b for any eigenfunction.

This finishes the proof of Lemma \ref{L:equiv}. \qed

\subsection{Equivalence $3\Leftrightarrow 4$}
Since conditions $3(a)$ and $4(a)$ are the same, it is sufficient
to prove equivalence of $3(b)$ and $4(b)$. Implication
$3(b)\Rightarrow 4(b)$ is straightforward. Indeed, one can choose
in (\ref{E:Darboux_eigenf}) instead of $\psi_\lambda$ one of the
eigenfunctions $\phi_{m,l}$ introduced in
(\ref{E:Laplace_eigenf_harm}), provided $\lambda\neq 0$ is a zero
of the Bessel function $J_{m+n/2-1}$. In this case, the integral
in (\ref{E:orthog}) evaluates to be proportional to
$\widehat{g}_{m,l}(\lambda)=\int\limits_S
\widehat{g}(\lambda,\theta)Y^m_l(\theta)d\theta$. Thus, vanishing
of these expressions for all $m,l$ and $\lambda$ as described, is
equivalent to the condition $4(b)$.

The converse implication $4(b)\Rightarrow 3(b)$ follows
analogously, if one takes into account the completeness of the
system of eigenfunctions $\phi_{m,l}$.

\subsection{Implication  $2+3+4\Rightarrow 1$}

Our goal here is, assuming any of the equivalent assumptions $2,
3, 4$ (or a combination of those), to show existence of a function
$f(x)$ supported in $B$ such that the restriction of its spherical
mean Radon transform $G(x,t)$ onto the lateral boundary $S\times
[0,2]$ of the cylinder $C$ coincides with the function $g$.

Using the Darboux equation reformulation that we have mentioned
before, this is equivalent to showing existence in $\RR^n \times
[0,\infty)$ of a solution $G(x,t)$ of the Darboux equation
(\ref{E:Darboux}) such that $G(x,0)=f(x)$, $G_t(x,0)=0$, and
$G\vert_{S\times [0,2]}=g$, for a function $f$ supported in $B$.
Then this function $f$ would be a pre-image under the spherical
mean transform $R$ of the data $g$.

Our strategy consists of solving the following sequence of
problems:
\begin{itemize}

\item Showing that the solution $G_+(x,t)$ of the interior problem
(\ref{E:Darboux}) is even and smooth with respect to $t$ on the
whole $t$-axis. This would, in particular, provide us with a
candidate $f(x)=G(x,0), x\in B$ for the pre-image.

\item Using rotational invariance, reducing the problem to single
spherical harmonic terms of $g$, $G$, and $f$.

\item Showing that each such term $G^m$ of $G$ extends to the
whole space $\RR^n$ as a global solution of Darboux equation.

\item Proving that the value $G^{m}(x,0)$ is supported inside the
ball $B$ and coincides with the corresponding harmonic term
$f_{m}$ of $f$. This will show that $Rf^{m}=g^{m}$.

\item Now an immediate continuity argument will show that $Rf=g$,
which will finish the proof of this implication, and thus of the
whole theorem.
\end{itemize}

Let us start realizing this program.

{\bf Interior solution $G_+$.}

Conditions 2(b) and 3(b) mean that the equivalent requirements
1(a) and 1(b) of Lemma \ref{L:equiv} are satisfied. Then the
second claim of this lemma guarantees that the interior solution
$G_+(x,t)$ can be continued to an even, smooth as $H^s(B)$-valued
function of $t$ solution of Darboux equation in the infinite
cylinder $B\times \RR$. This resolves the first step of our
program.

{\bf Projection to the $O(n)$- irreducible representations.}

Each irreducible sub-representation $X^m$ of the representation of
the orthogonal group $O(n)$ on functions on $\RR^n$ by rotations
consists of homogeneous harmonic polynomials of a fixed degree $m$
(e.g., \cite{Stein,Vilenkin}). Restrictions of the elements of
$X^m$ to the unit sphere $S$ are spherical harmonics of degree
$m$. We denote, as before, by $Y^{m}_{l}(\theta), l
=1,\cdots,d(m), m=1, \cdots$ an orthonormal basis in $X^m$. The
orthogonal projection $L^2(S) \mapsto X^m$ will be denoted by
$\mathcal{P}^m$:
$$
(\mathcal P^m h)(x)=\int_S h(y)Z^m_x(y)dS(y),
$$
where $Z^m_x$ is the zonal spherical harmonic of degree $m$ with
the pole $x$ (e.g., \cite[Chapter 4, Section 2]{Stein}).
Since Bessel operator $\mathcal{B}$ and
Laplace operator $\Delta=\Delta_x $ both commute with the action
of $O(n)$, the projection onto $X^m$
$$
G^m(x,t)= \left(\mathcal{P}^m G\right)(x,t)
$$
reduces the Darboux equation. I.e., $G^m$ solves the same Darboux
equation with the zero data for $t=2$ and with the boundary data
$g^{m}=\mathcal P^{m}g$. Clearly, we also have $G^m(x,0)=f^m(x)$.
So, let us assume for now that $G=G^{m}$, $g=g^{m}$, and $f=f^m$.

Since functions $Y^m_l,l=1,\cdots,d(m)$ form an orthonormal basis
of $X^m$, we have
\begin{equation}\label{E:special_form}
\begin{array}{c}
G^m(x,t)=\sum\limits_l^{d(m)} g_l(r,t) Y_l^{m}(\theta),\\
g^m(\theta,t)=\sum\limits_l^{d(m)} g_l(t) Y_l^{m}(\theta).
\end{array}
\end{equation}

As we have already seen, the Fourier-Bessel transform takes the
solution $G^m$ of Darboux equation to a function $\widehat
G^m(x,\lambda)$ that satisfies the equation
$$
\Delta_x \widehat G^m(x,\lambda)=-\lambda^2 \widehat
G^m(x,\lambda)
$$
in $B$. Due to this and the special form (\ref{E:special_form}) of
$G^m$, its Fourier-Bessel transform can be written as
\begin{equation}\label{E:FB}
\widehat G^m(r\theta,\lambda)= j_{n/2-1+m}(\lambda r)(\lambda r)^m
\sum\limits_{l=1}^{d(m)} b_l(\lambda)Y_l^m(\theta),
\end{equation}
and correspondingly
\begin{equation}\label{E:FB_boundary}
\widehat g^m(\theta,\lambda)= j_{n/2-1+m}(\lambda )
(\lambda)^m\sum\limits_{l=1}^{d(m)} b_l(\lambda)Y_l^m(\theta).
\end{equation}
Now observe that the right hand side of (\ref{E:FB}) is defined
for all $r$, not only for $r \le 1$ and therefore defines a smooth
extension of $\widehat G^m(r\theta,\lambda)$  for $r>1.$ Thus, we
can think of $\widehat G^m(x,\lambda)$ as smooth in $x=r\theta$
function defined for {\bf all} $(x,\lambda) \in \RR^n\times \RR$.
For $x \in B$, due to Lemma \ref{L:PW-Bessel}, this function is of
the Paley-Wiener class in $\lambda$, being Fourier-Bessel
transform in $t$ of the compactly supported smooth function
$G^m(x,t)$. However, at this stage we do not know much about its
behavior with respect to $\lambda$ for $x \notin B$. To gain this
knowledge, we need some control over the smoothness and growth of
the coefficients $b_l(\lambda)$, which are defined for all
$\lambda \in \CC$.

Computing the Fourier coefficients with respect to the orthonormal
basis $Y^m_l$ of spherical harmonics, we obtain
\begin{equation}\label{E:lambda}
b_l(\lambda)(\lambda r)^m j_{n/2-1+m}(\lambda
r)=\int\limits_{\theta \in S} \widehat G^m(r\theta,
\lambda)Y_l^m(\theta) dS(\theta).
\end{equation}
The functions $b_l(\lambda)$ are clearly analytic at any point
$\lambda_0 \neq 0$. Indeed, for any such $\lambda_0$ one can
choose $r<1$ so that $j_{n/2-1+m}(\lambda r) \neq 0$ for $\lambda
$ near $\lambda_0$. Thus, $b_l$ is analytic near $\lambda_0$ as
the ratio of two analytic functions with non-vanishing
denominator.

This argument does not work at $\lambda=0$. Moreover, smoothness
of $b_l$ at $\lambda=0$ is not guaranteed immediately by
(\ref{E:lambda}), and requires the moment condition. Indeed, let
us restrict (\ref{E:lambda}) to the boundary $S$ to get
\begin{equation}\label{E:lambda_boundary}
b_l(\lambda)(\lambda)^m j_{n/2-1+m}(\lambda )=\int\limits_{\theta
\in S} \widehat g^m(\theta, \lambda)Y_l^m(\theta) dS(\theta).
\end{equation}
As we already established in Lemma \ref{L:moments_interpretation},
the moment condition is equivalent to the integral in the right
hand side in (\ref{E:lambda_boundary}) vanishing at $\lambda=0$ to
the order at least $m$. On the other hand, the function
$\lambda^mj_{n/2-1+m}(\lambda)$ in (\ref{E:lambda_boundary}) has
zero of order $m$ at $\lambda=0$. Then, dividing by $\lambda^m
j_{n/2-1+m}(\lambda)$
 in (\ref{E:lambda_boundary}), we conclude that $b_l(\lambda)$ is
smooth at $\lambda=0$.

Thus, $b_l$ is an entire function. Now we need to estimate its
growth at infinity (looking for Paley-Wiener estimates). Due to
the Paley-Wiener class estimates that we have for the expression
in the right hand side of (\ref{E:lambda_boundary}) and known
behavior of Bessel functions, it is a standard exercise to show
that their ratio $b_l(\lambda)$ is of the same Paley-Wiener class
as the numerator. Indeed, this was treated in
\cite{AmbKuch_range}. The estimate from below for Bessel functions
provided in \cite[Lemma 6]{AmbKuch_range} and its consequent usage
there show that this in fact is true\footnote{An alternative proof
can be found in Section \ref{S:lemmas}}.

{\bf Extending $G^m$}

Now we can apply the inverse Fourier-Bessel transform in $\lambda$
to the extended function $\widehat G^m(x,\lambda), x \in \RR^n$ in
(\ref{E:FB}). We use the same notation $G^m(x,t)$ for the obtained
function. This is justified, since by the construction this
function satisfies the Darboux equation and it coincides with the
original interior ``$m$-irreducible'' solution $G^m(x,t)$ in the
cylinder $B \times \RR$.  One can also observe that, due to the
Paley-Wiener Lemma \ref{L:PW-Bessel}, it is smooth with respect to
$t$ and supported in $B\times [-2,2]$.

{\bf The size of the support of $G^m(x,0)$.}

Using the relation $G^m(x,0)=f^m(x), x\in B$ and applying the
inverse Fourier-Bessel transform to (\ref{E:FB}), one finds
\begin{equation}\label{E:f_as IFT}
f^m(x)=f^m(r\theta)={\rm const}\int\limits^\infty_0 (\lambda r)^m
j_{n/2-1+m}(\lambda r)\left(\sum\limits_{l=1}^{d(m)}b_l(\lambda)
Y^m_l(\theta)\right)\lambda^{n-1}d\lambda .
\end{equation}
We notice that now we can define an extension $F^m(x)$ of $f^m(x)$
to the whole space $\RR^n$ by applying the formula (\ref{E:f_as
IFT}) for $r>1$:
\begin{equation}\label{E:F}
F^m(x)=F^m(r\theta):={\rm const}\int\limits^\infty_0 (\lambda r)^m
j_{n/2-1+m}(\lambda r)\left(\sum\limits_{l=1}^{d(m)}b_l(\lambda)
Y^m_l(\theta)\right)\lambda^{n-1}d\lambda,
\end{equation}
which, according to known results (e.g., \cite[Ch. IV, Theorem
3.10]{Stein}), is just the inverse $n$-dimensional Fourier
transform of the following function
\begin{equation}\label{E:H}
H(x)=H(\lambda
\theta)=\lambda^m\sum\limits_{l=1}^{d(m)}b_l(\lambda)
Y^m_l(\theta)
\end{equation}
 in $\RR^n$, written in polar coordinates $x=\lambda\theta$.
Consider for a moment only real values of $\lambda$. One sees
immediately that, due to the Paley-Wiener estimates on $b_l$,
function $H(x)$ is smooth outside the origin and decays with all
its derivatives faster than any power of $|x|=|\lambda|$. If we
show smoothness at the origin, then $H$ will be proven to belong
to the Schwartz class. According to standard considerations of
Radon transform theory \cite{GGG1, GGG2, Helg_Radon}, this
function is smooth at the origin if and only if for any $k$ the
expression
$$
\frac{\partial ^k H(\lambda \theta)}{\partial \lambda
^k}|_{\lambda=0}=\sum\limits_{l=1}^{d(m)}\frac{\partial ^k
(\lambda^m b_l(\lambda))}{\partial \lambda
^k}|_{\lambda=0}Y^m_l(\theta)
$$
is the restriction to the unit sphere of a homogeneous polynomial
of degree $k$ with respect to $\theta$. Due to the form of the
last expression, this means that $\frac{\partial ^k (\lambda^m
b_l(\lambda))}{\partial \lambda ^k}(0)=0$ for any $k$ such that
either $k<m$ or $k-m$ is odd. The case $k<m$ is obvious, due to
smoothness of $b_l$ and the presence of the factor $\lambda^m$.
Due to the structure of the functions $\hat g (\lambda)$ and
Bessel functions, discussed already, the condition for $k-m$ odd
is automatic.

Hence, $H$ belongs to the Schwartz space. Then its inverse Fourier
transform, which we previously denoted by $F^m(x)$, is in Schwartz
space itself. We now need to establish that $F^m$ is supported
inside $B$.

By its construction, $F^m$ is the value at $t=0$ of a global
solution $G^m$ of the Darboux equation. Since $G^m|_{S\times
\RR^+}=g^m$, Asgeirsson Theorem \cite{Asgeirsson,John} implies
that $RF^m=g^m$. In particular, the integrals of $F^m$ over all
spheres centered inside $B$ and of radii  $t\geq 2$ are equal to
zero. Indeed, such an integral over a sphere centered at $x\in B$
of radius $t$ is equal to $G^m(x,t)$, which is known to be zero by
construction of $G^m$.

Lemma 2.7 in \cite[Ch. 1]{Helg_Radon} claims that if a function
decays faster than any power and its integrals over all spheres
surrounding a convex body $B$ are equal to zero, then the function
is zero outside of $B$. In our case we do not have all such
spheres, but only the ones of radii at least $2$ and centered in
$B$, where $B$ is the unit ball. However, a simple exercise is to
check that the proof of the cited lemma still holds and thus $F^m$
is supported in $B$. This means that in fact $F^m(x)$ is the zero
extension of $f^m (x)$ outside the ball $B$.

{\bf Final step: proving $Rf=g$}

We already have constructed a function $f\in C^\infty_0(B)$ such
that for each its component $f^m$ corresponding to an irreducible
representation $X^m$, the equality $Rf^m=g^m$ holds. Since we have
expansions $f=\sum_m f^m$ and $g=\sum_m g^m$ converging in any
Sobolev space, the equality $Rf=g$ immediately follows by
continuity.

This finishes the proof of implication (2)+(3)+(4) $\Rightarrow$
(1) and thus completes the proof of Theorem \ref{T:Main_ball}.
\qed

\section{Proof of Theorem \ref{T:Main_odd}}\label{S:odd}
We assume now that the dimension $n$ is odd. The claim of Theorem
\ref{T:Main_odd} is that the statement of Theorem
\ref{T:Main_ball} can be proven without using moment conditions of
Lemma \ref{L:moment_ball} or Lemma \ref{L:moment_general}. In
other words, we claim that in odd dimensions, the moments
conditions (2a)-(4a) of Theorem \ref{T:Main_ball} follow from the
equivalent orthogonality conditions (2b)-(4b). Of course, this
effects only the sufficiency part of this theorem, since we have
proven that in any dimension the moment conditions are necessary
for $g$ to be in the range. Thus, in odd dimensions the moment
conditions are redundant.

In order to prove Theorem \ref{T:Main_odd}, we will follow the
part of the proof of Theorem \ref{T:Main_ball} where the moment
conditions were not used, and then will finish the proof avoiding
the moment conditions.

As we have proven in Section \ref{S:proof_ball}, the equivalent
conditions (2) - (4) of Theorem \ref{T:Main_odd} (which coincide with
conditions (2b) - (4b) of Theorem \ref{T:Main_ball}), imply the
existence of a smooth solution $G_+(x,t)$ of Darboux equation in
the solid cylinder $C=B \times \RR$, even with respect to $t$ and
such that $G_+(x,t)=0$ for $|t|\geq 2$ and $G_+|_{S\times\RR}=g$.
Here $g=g(x,t)$ is the function introduced in Theorem
\ref{T:Main_odd}. We emphasize again that the moment conditions
were not used in this derivation.

The key point in the proof of Theorem \ref{T:Main_odd} is the
following auxiliary statement, which in odd dimensions can be
proven without using moment conditions:

\begin{proposition}\label{P:smooth}
The function $G_+$ and all its partial derivatives in $x$ vanish
on the sphere $S_0= \{|x|=1, t=0\}\subset \RR^n \times \{0\}.$
\end{proposition}

We will postpone the proof of this Proposition, and will show now
that it implies the statement of Theorem \ref{T:Main_odd}.

\subsection{Derivation of Theorem \ref{T:Main_odd} from Proposition
\ref{P:smooth}}\label{SS:global}

Let $f(x):=G_+(x,0)$. This functions is smooth in $B$ and,
according to Proposition \ref{P:smooth}, vanishes to the infinite
order at $S=\partial B$. Hence, the function obtained by zero
extension of $f(x)$ to the whole $\RR^n$, is smooth. We will use
the same notation $f(x)$ for this function:
$$
f(x):=\begin{cases} f(x) \mbox{ for } |x|\leq 1\\ 0 \mbox{ for }
|x|\geq 1.\end{cases}
$$
Let us now define
$$
G(x,t)=\frac{1}{\omega}\int\limits_{y \in S}f(x+ty)dS(y),\,x\in
\RR^n
$$
as the spherical mean transform of this function $f$. Theorem
\ref{T:Main_odd} will be proven, if we show that $G=G_+$ inside
$C$. Indeed, then $Rf|_{S\times [0,\infty)}=g$, and thus $g$ is in
the range.

The function $G$ solves Darboux equation in the entire space, and
therefore the difference
$$
Z(x,t)=G_+(x,t)-G(x,t)
$$
solves this equation in the cylinder $B \times [0,\infty)$ (in
fact, it extends to an even solution in $B \times \RR$).

\begin{figure}[ht]
\begin{center}
\scalebox{0.5}{\includegraphics{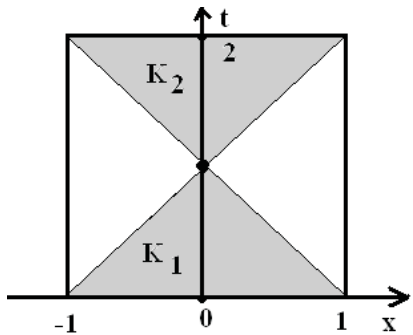}}\\
\caption{The cones $K_1$ and $K_2$.}
\label{F:cones}
\end{center}
\end{figure}

Our goal is to prove that the solutions $G$ and $G_+$ agree on the
cylinder $C$. At the initial moment $t=0$, we have for $x \in B$:
$$
\begin{array}{c}
Z(x,0)=f(x)-f(x)=0,\\
Z_t(x,0)=G_{+,t}(x,0)-G_t(x,0)=0.
\end{array}
$$
Then $Z$ vanishes inside the corresponding characteristic cone:
$$
Z(x,t)=0, (x,t) \in K_1:=\{0 \le t \le 1-|x| \}.
$$
On the other hand, the interior solution $G_+(x,t)$, according to
its construction, vanishes for $t \ge 2$. The spherical means
$G(x,t)$ also vanish for $|x|\leq 1, t\geq 2$, since then the
sphere $\{y: |y-x|=t \}$ does not intersect $B$, and hence also
the support of $f$. Thus, the difference $Z=G_+-G$ vanishes in the
cylinder $\{|x|\leq 1, t\geq 2\}$. By the same dependence domain
argument, $Z$ vanishes on the backward characteristic cone $K_2$:
$$
Z(x,t)=0, (x,t) \in K_2:=\{|x|+1 \le t \le 2 \}.
$$
Thus, the difference $Z$ of the two solutions vanishes on the
union of two characteristic cones:
$$
Z(x,t)=0, (x,t) \in K_1 \cup K_2
$$
that have the common vertex $(0,1)$ (see Figure \ref{F:cones}).

Notice that the union $K_1 \cup K_2$ of the two cones contains the
segment $x=0, 0 \le t \le 2$. Each point $(0,t)$ of this segment,
except the vertex $(0,1)$, belongs to the interior of $K=K_1 \cup
K_2$. Since $Z$ is smooth and $Z\vert_K=0$, all partial
derivatives $\partial_x^{\alpha}Z(0,t)$ vanish for all $t \neq 1$.
By smoothness, this is also true for the vertex $(0,1)$. Thus, function $Z$
 vanishes, along with all its derivatives, on the entire line $x=0$. We claim
that this implies that $Z =0$ for all $(x,t) \in C$. The following
lemma does the job:
\begin{lemma}\label{L:cones}
Let $Z(x,t)$ be an even and compactly supported in $t$ smooth
solution of Darboux equation in $B \times \RR$ and $x_0$ an
interior point in $B$. If $(\partial_x^{\alpha}Z)(x_0,t)=0$ for
any multi-index $\alpha$ and $t \in \RR,$ then $Z(x,t)=0$ for all
$(x,t) \in B \times \RR.$
\end{lemma}
\pf Let us apply Fourier-Bessel transform with respect to $t$ to
the function $Z(x,t)$. Then Darboux equation transforms to
Helmholtz equation and the resulting function $\widehat
Z(x,\lambda)$ is an eigenfunction of Laplace operator:
$$
\Delta \widehat Z(x,\lambda)=-\lambda^2 \widehat Z(x,\lambda).
$$
By the condition of the lemma, $\widehat Z(x,\lambda)$ has at
$x_0$ a zero of infinite order. Since eigenfunctions of the
Laplace operator $\Delta$ are known to be real-analytic, we
conclude that $\widehat Z(\cdot,\lambda)\equiv 0$. Taking inverse
Fourier-Bessel transform, we get $Z\equiv 0$. \qed

This finishes the proof of Theorem \ref{T:Main_odd}, modulo the
proof of Proposition \ref{P:smooth}, which we provide in the next
sub-section.

\subsection{Proof of Proposition \ref{P:smooth}}

The goal of this subsection is to prove that, as the Proposition
\ref{P:smooth} states, at the boundary points $(x,0)$, the
infinite order zero of the boundary data $g(x,t)$ at $t=0$ implies
infinite order zero of $G(x,t)$ with respect to $x$ at
$t=0,|x|=1$.

This proof will be close to the proof of Proposition 7 in
\cite{FinRak}. Following \cite{FinRak}, we will reduce the
initial-boundary value problem for Darboux equation to a problem
for the one dimensional wave equation. Again, as in \cite{FinRak},
we will use separation of variables in polar coordinates and
dimension reduction. However, we will do this in a somewhat
different manner. Besides, since we are dealing with Darboux
equation rather than with the wave equation, an additional
integral transform will be required with respect to $t$.

We will break the proof into several steps.

\subsubsection{Separation of variables.}
First of all, we decompose the solution $G$ into spherical
harmonic parts $G^m$ that belong to the irreducible
representations $X^m$ of the rotation group $O(n)$:
$$
G(x,t)=\sum\limits_m
G^m(x,t)=\sum\limits_m\sum\limits_{l=0}^{d(m)} \Psi_l(r,t)r^m
Y_l^m(\theta), x=r\theta.
$$
It is clear that it is sufficient to prove the claim of the
Proposition for each of these components
$$
G^m(x,t)=\sum\limits_{l=0}^{d(m)} \Psi_l(r,t)r^m Y_l^m(\theta),
x=r\theta.
$$
Indeed, we know that $G$ belongs to any Sobolev space $H^s(B)$ and
hence the convergence of the spherical harmonics expansion and
Sobolev embedding theorems will deduce the claim for $G$ from
those for each $G^m$.

So, we will assume from now on that
$G=G^m=\sum\limits_{l=0}^{d(m)} \Psi_l(r,t)r^m Y_l^m(\theta)$.
Correspondingly, we also expand the data $g$ and assume that
$$
g(\theta,t)=\sum\limits_{l=0}^{d(m)} \Psi_l(t) Y_l^m(\theta).
$$
The Darboux equation for $G(x,t)$ implies that the coefficients
$\Psi(r,t)=\Psi_l(r,t)$ satisfy the following PDE:
\begin{equation}\label{E:Darboux_r,t}
\partial_t^2 \Psi +\frac{n-1}{t} \partial_t \Psi=
\partial_r^2 \Psi + \frac{n-1+2m}{r} \partial_r \Psi.
\end{equation}
In order to prove the Proposition, it suffices to prove that
$$
\frac{\partial^p \Psi}{\partial r^p} (1,0)=0, p=0,1,\cdots.
$$

So, we will concentrate now on studying the solution $\Psi$ of
(\ref{E:Darboux_r,t}). By its construction, $\Psi(r,t)$ satisfies
the following conditions:

\noindent a) $\Psi(r,t)$ is smooth end even with respect to $t \in
\mathbb R$ and $r \in [-1,1]$,

\noindent b) $\Psi(r,t)=0$ for $|t|\geq 2$,

\noindent c) $\Psi(1,t)$ vanishes at $t=0$ to the infinite order.

\subsubsection{Reduction to the one dimensional wave equation.}
Let us set $p=\frac{n-2+2m}{2}$ and apply the inverse Poisson
transform $\mathcal{P}_p^{-1}$ (see (\ref{E:Poisson_inversion}))
with respect to the variable $r$:

$$
\left(\mathcal{P}^{-1}\Psi\right)(r,t)=(\mbox{const})
r\left(\frac{\partial}{\partial
(r^2)}\right)^{(n-1+2m)/2}r^{n-2+2m} \Psi(r,t).
$$
This will reduce the right hand side of equation
(\ref{E:Darboux_r,t}) to the second $r$-derivative. We also apply
Weyl transform $\mathcal{W}_{(n-2)/2}$ (\ref{E:Weyl}) in the
variable $t$. As the result, we obtain the function
$$
U(r,t)=\left(\mathcal{P}_{(n-2+2m)/2,r}^{-1}\mathcal{W}_{(n-2)/2,t}\Psi\right)(r,t),
$$
which, due to the intertwining properties (\ref{E:intertwine}), solves the one dimensional wave equation
\begin{equation}\label{E:wave2}
U_{tt}-U_{rr}=0.
\end{equation}
We can observe now that this new function $U$ has the following properties:

\noindent 1) $U(\pm r,\pm t)=U(r,t),$

\noindent 2)$U(r,t)=0$ for $|t| \ge 2.$

\noindent 3)
$U(1,t)=q(t):=\left(\mathcal{P}^{-1}\mathcal{W}\Psi\right)(1,t)$
and $q(t)=0, |t| \ge 2.$

The evenness property 1) follows from
the evenness of $\Psi$ and the fact that both transforms we applied preserve it. Property 2) follows from the analogous property b) of $\Psi$ and the properties of
the Weyl transform $\mathcal{W}$. Now, 3) is just the definition of the
boundary values of $U$ for $r=\pm 1$ combined with 2).

The unique solution of (\ref{E:wave2}) in the domain $-1 \le r \le
1, t \geq 0$ with these properties, is the sum of the following
two progressing waves:
\begin{equation}\label{E:progressing}
U(r,t)=q(t-r+1)+q(t+r+1).
\end{equation}
Indeed, this sum clearly satisfies (\ref{E:wave2}). If $r=\pm 1, t
\ge 0,$ then $U(\pm 1,t)=q(t)+q(t+2)=q(t)$,  due to 3). So,
property 2) holds. Also, $q(t-r+1)+q(t+r+1)$ is obviously even in
$r$. Taking into account that $q(t)=0$ for $|t|>2$, implies that
$q(t-r+1)+q(t+r+1)=0$ for $t>2$. Thus, by standard uniqueness
theorem, the solutions $U(r,t)$ and $q(t-r+1)+q(t+r+1)$ coincide.

The representation (\ref{E:progressing}) has important consequences. The first is given in the following
\begin{lemma}\label{L:U^p(0)=0}
For all $p \in \mathbb Z_+$, one has
$$
\partial_r^pU(\pm 1,0)=0.
$$
\end{lemma}
\pf Let us show first that  the evenness of $U$ with respect to $t$
implies infinite order zero at 0 of the boundary value $q$. Indeed, observe that
for any odd number $p=2s-1$, one has
\begin{equation}\label{E:q}
q^{(2s-1)}(-r+1)+q^{(2s-1)}(r+1)=(\partial_t^{2s-1}U)(r,0)=0.
\end{equation}
Substituting  $r=1$ along with using 3) leads to
$$
q^{(2s-1)}(0)= q^{(2s-1)}(0)+q^{(2s-1)}(2)=0.
$$
On the other hand, differentiation of (\ref{E:q}) with respect to $r$ implies
$$
-q^{(2s)}(-r+1)+q^{(2s)}(r+1)=0.
$$
Again, substituting $r=1$ yields
$$
-q^{(2s)}(0)=-q^{(2s)}(0)+q^{(2s)}(2)=0.
$$
Thus, all derivatives $q^{(p)}(0)$ vanish. This, in turn, implies that
$U(r,t)$ has a zero of infinite order (with respect to $r$) at $(1,0)$:
\begin{equation}\label{E:zerosU}
\partial_r^p U(1,0)= (-1)^p q^{(p)}(0)+q^{(p)}(2)=0.
\end{equation}
Since $U$ is even in $r$, also $\partial_r^pU(-1,0)=0.$ \qed

\subsubsection{Proving that $\Psi$ vanishes to infinite order at $|x|=1, t=0$.}

We can now finish the proof of the Proposition \ref{P:smooth} by
showing that $\partial_r^p \Psi(1,0)=0$ for any $p=0,1,\cdots$.
Notice, that it suffices to check this identity only for even $p$.
Indeed, $\Psi(r,t)$ vanishes at $(1,0)$ with all derivatives with
respect to $t$. Then the equation (\ref{E:Darboux_r,t}) implies
that all iterates of the Darboux operator acting in the variable
$r$ vanish:
$$
(\partial_r^2 +((n-1+2m)/2)\partial_r)^N\Psi(1,0)=0.
$$
Now vanishing of all even order derivatives
$\partial_r^{2k}\Psi(1,0)$ clearly implies vanishing  of the
derivatives of odd order as well.

It will be convenient to use the following simple relation:
\begin{lemma}\label{L:v}  For any smooth even function $v(t)$ the following relation
holds:
$$
\left(\frac{d}{d(t^2)}\right)^j v(0)=\left(\frac{1}{2t}
\frac{d}{dt}\right)^j v(0)=\frac{j!}{(2j)!}v^{(2j)}(0).
$$
\end{lemma}
\pf This equality follows from the Taylor formula. \qed

Let us write now the relation
$\Psi=(\mathcal{P}\mathcal{W}^{-1})U$, using (\ref{E:Poisson}) and
(\ref{E:Weyl}), explicitly:
\begin{equation}\label{E:Psi-U}
\Psi(r,t)=\mbox{const}\int\limits_{-1}^1
\left(\frac{\partial}{\partial (t^2)}\right)^{(n-1)/2} U(\mu
r,t)(1-\mu^2)^{(n-3)/2}d\mu.
\end{equation}
According to the property c) of the function $\Psi$, its all
$t$-derivatives at $(1,0)$ are equal to zero. The odd order
derivatives vanish due to the evenness of $\Psi$, so only the even
order derivatives carry interesting information for us:
$$
\partial_t^{2j}\Psi(1,0)=0, p=0,1,\cdots .
$$
Let us translate, using (\ref{E:Psi-U}), these equalities into the
language of function $U$. Differentiating $2j$ times with respect
to $t$ at the point $r=1,t=0$ under the sign of the integral in
(\ref{E:Psi-U}) leads to the expression
$\partial_t^{2j}\partial_{t^2}^{(n-1)/2} U(\mu ,0)$. Since $U(\mu
,t)$ is even with respect to $t$, Lemma \ref{L:v} gives
$$
\partial_t^{2j}\partial_{t^2}^{(n-1)/2} U(\mu ,0)=\mbox{const}
\partial_t^{2j+n-1} U(\mu,0).
$$
Taking into account that $2j+n-1$ is an even number, the wave
equation (\ref{E:wave2}) (or the progressing wave expansion
(\ref{E:progressing})) yields
\begin{equation}\label{E:t_r}
\partial_t^{2j+n-1}U(r,0)=\partial_r^{2j+n-1}U(r,0).
\end{equation}

Thus,
\begin{equation}\label{E:U_mu}
\partial_t^{2j}\Psi(1,0)=\mbox{const} \int\limits_{-1}^1 \partial_{\mu}^{2j+n-1} U(\mu,0)(1-\mu^2)^{(n-3)/2}d\mu=0,
\end{equation}
for all $j=0,1,\cdots.$

We can now reformulate the identities (\ref{E:U_mu}) as follows:
\begin{lemma}\label{E:B} The function $B(\mu)=(\partial_{\mu}^{(n+1)/2}U)(\mu,0)$ is orthogonal
in $L^2[-1,1]$ to all polynomials $P$ of degree $\deg P \le
(n-3)/2$ of the same parity as the natural number $(n-3)/2$ (i.e.
$P$ is even if $(n-3)/2$ is even, and odd if $(n-3)/2$ is odd).
\end{lemma}
\pf Integration by parts $2j+(n-3)/2$ times in (\ref{E:U_mu}) and
vanishing of derivatives of $U$ at $(1,0)$ lead to
\begin{equation}\label{E:Legendre}
\int\limits_{-1}^1 B(\mu) L_{\frac{n-3}{2}}^{(2j)}(\mu)d\mu=0,
\end{equation}
where $L_m(\mu)=\mbox{const} \partial_{\mu}^m (\mu^2-1)^m$ are
Legendre polynomials. Since the derivatives $L_m^{(2j)},
j=0,1,\cdots $ clearly span the space of polynomials $P$ of degree
$\deg P \le m$ that have the same parity as $m$, this proves the
Lemma. \qed

Let us now prove that $r$-derivatives of $\Psi$ vanish at $(1,0)$.
We differentiate $2k$ times the identity (\ref{E:Psi-U}) in the
variable $r$. The wave equation (\ref{E:wave2}) and Lemma \ref{L:v}
imply
$$
(\partial_{\mu}^{2k}\partial_{t^2}^{(n-1)/2}U)(\mu,0)=
\mbox{const}(\partial_t^{2k+n-1}U)(\mu,0).
$$
Then (\ref{E:Psi-U}) leads to
$$
\begin{array}{c}
\partial_r^{2k}\Psi(1,0)=\mbox{const}\int\limits_{-1}^1 (\partial_{\mu}^{2k} \partial_{t^2}^{(n-1)/2}U)(\mu,0)\mu^{2k}
(1-\mu^2)^{(n-3)/2}d\mu\\
=\mbox{const}\int\limits_{-1}^1(\partial_{\mu}^{2k+n-1}U)(\mu,0)\mu^{2k}(1-\mu^2)^{(n-3)/2}d\mu.
\end{array}
$$
Integration by parts $2k+(n-3)/2$ times gives
$$
\partial_r^{2k}\Psi(1,0)=\mbox{const}\int\limits_{-1}^1(\partial_{\mu}^{(n+1)/2}U)(\mu,0) P(\mu)d\mu=
\mbox{const} \int\limits_{-1}^1 B(\mu)P(\mu)d\mu,
$$
where
$$P(\mu)=\partial_{\mu}^{2k+(n-3)/2}(\mu^{2k}(1-\mu^2)^{(n-3)/2})$$
is the polynomial of degree $(n-3)/2$ and of the same parity as
$(n-3)/2$. Now Lemma \ref{E:B} claims that all such integrals are
equal to zero. Thus, $\partial_r^p\Psi(\pm 1,0)=0$, which finishes
the proof of Proposition \ref{P:smooth}. \qed

\section{Proof of Theorem \ref{T:estimates}}

The necessity was already established in (\ref{E:Q_estimate}). To
prove sufficiency, let us check that the conditions 3(a) and 3(b)
of Theorem \ref{T:Main_ball} hold, which will imply that $g$ is in
the range of the transform $R$.

First of all, the moment condition 3(a) is obviously weaker than
our condition.

To check 3(b), observe that the Fourier-Bessel transform
$$
\hat g(x,\lambda)=\int\limits_0^{\infty}g(x,t)j_{n/2-1}(\lambda t)
t^{n-1}dt,
$$
which is an entire function of $\lambda$, expands, according to
(\ref{E:Bessel_expansion}), into the power series with respect to
the spectral parameter $\lambda$:
$$
\hat g(x,\lambda)=\int\limits_0^{\infty}g(x,t)j_{n/2-1}(\lambda
t)t^{n-1}dt= \sum\limits_{k=0}^{\infty}C_k \lambda^{2k}M_k(x).
$$
If we replace here the functions $M_k(x)$ by their  extensions
$Q_k(x)$ in the unit ball $B$, then, due to estimates
(\ref{E:ball_estimate}) for $Q_k$, the extended series uniformly
converges in $B$ to a function
$$
\psi_{\lambda}(x):=\sum\limits_0^{\infty} C_k \lambda^{2k}Q_k(x).
$$
This function $\psi_{\lambda}$ is an eigenfunction of Laplace
operator:
$$
\Delta \psi_{\lambda} =\sum\limits_{k=0}
 ^{\infty}C_k \lambda^{2k} c_k Q_{k-1}=
 \sum\limits_{s=0}^{\infty} C_{s+1}c_{s+1}\lambda^{2(s+1)}Q_s=-\lambda^2 \psi_{\lambda},$$
due to the relation
$$
C_{s+1}c_{s+1}=-C_sc_s.
$$
This relation between the coefficients can be easily verified
using their explicit values
$$
C_k=(-1)^k\frac{\Gamma(p+1)}{2^{2k}k!\Gamma(p +k+1)},
c_k=2k(2k+n-2), p=\frac{n-2}{2}.
$$
Therefore, the function $\widehat g_{\lambda}(x)= \widehat
g(x,\lambda)$ extends to the ball $B$ as an eigenfunction of
Laplace operator with the eigenvalue $-\lambda^2$. Let us now
apply Stokes formula to $\psi_\lambda$ and a Dirichlet
eigenfunction $\phi=\phi_{\lambda}$. Taking into account that
$\hat g_{\lambda}=\psi_{\lambda}$ on $S$, we obtain
$$
\int\limits_S  \hat g_{\lambda} \partial_{\nu} \phi_{\lambda}dS=
  \int\limits_B (\psi_{\lambda}\Delta \phi_{\lambda} - \phi_{\lambda}\Delta \psi_{\lambda})dV=0.
$$
This provides the orthogonality condition 3(b)
(formula(\ref{E:orthog})). Thus, according to Theorem 10, function
$g$ belongs to the range of transform $R$. \qed

\section{Proof of Theorem \ref{T:sobolev}}
First of all, the necessity part of the proof of Theorem
\ref{T:Main_ball} clearly survives in the Sobolev case. Thus, we
only need to establish that a function $g\in
H_{s+(n-1)/2}^{\mbox{comp}}(S\times(0,2))$ satisfying the range
conditions, does belong to the range of $R_S$ on
$H_s^{\mbox{comp}}(B)$. Here the following stability estimate is
essential:
\begin{proposition}\label{P:stability}
For any $\varepsilon >0$, there exists a constant $C_\varepsilon$
such that for any $f\in H_s^0(B_{1-\varepsilon})$ the following
estimate holds:
\begin{equation}\label{E:stability}
   C^{-1}_\varepsilon \|R_S
    f\|_{H_{s+(n-1)/2}^0(S\times(\varepsilon,2-\varepsilon))}\leq \|f\|_{H_s^0(B_{1-\varepsilon})}\leq C_\varepsilon \|R_S
    f\|_{H_{s+(n-1)/2}^0(S\times(\varepsilon,2-\varepsilon))}.
\end{equation}
Here $B_{1-\varepsilon}$ is the ball of radius $1-\varepsilon$
centered at the origin.
\end{proposition}
This proposition combines the known injectivity of the operator
$R_S$ on compactly supported functions (see, e.g., \cite{AQ}) with
the elliptic estimates obtained recently in \cite{Pal_funk}. Such
estimates can also be derived from the known ellipticity theorem
for the pseudo-differential normal operator $R^*_S R_S$ and FIO
results of \cite{Hor}. The ellipticity theorem was obtained in
\cite{GS} (see also \cite{GU}--\cite{Guill_Ster}, \cite{Q1980}).

Let us now have a function $g(x,t)\in
H_{s+(n-1)/2}^{\mbox{comp}}(S\times(0,2))$ satisfying any of the
range conditions formulated in the theorem. We can choose a
positive $\varepsilon$ such that $g\in
H_{s+(n-1)/2}^0(S\times(2\varepsilon,2-2\varepsilon))$. Let us
extend $g$ to a function $g_1(x,y)$ on $\RR^n\times\RR^n$, radial
with respect to $y$, as follows:
$$
g_1(x,y)=g(x,|y|).
$$
The proof of Theorem \ref{T:Main_ball} and lemmas preceding it
show that the range conditions require the following: if one
expands $g_1$ into a series of spherical harmonic terms with
respect to $x$, and then takes the $n$-dimensional Fourier
transform with respect to $y$ of each term, the resulting
functions have zeros at certain prescribed locations and of
prescribed multiplicities. Let us now consider an even, smooth,
compactly supported radial approximation $\psi_k$ of the
delta-function in $\RR^n_y$, such that $\psi_k \rightarrow \delta$
in distributional sense when $k \rightarrow \infty$ and that the
support of $\psi_k$ shrinks to $\{0\}$ when $k \rightarrow
\infty$. Then the convolution $g_k(x,y)=\psi_k (y)*_y g_1(x,y)$
for a large $k$ is a $C_0^\infty$, radial with respect to $y$
function. If we now denote $g_k(x,t)=g_k(x,y)$ for $|y|=t$, we get
for large $k$ a smooth function with support in $S\times
[\varepsilon,2-\varepsilon]$ and such that $g_k \rightarrow g$ in
$g\in H_{s+(n-1)/2}^0(S\times(\varepsilon,2-\varepsilon))$. Since
the Fourier transform of the convolution is the product of Fourier
transforms, we see that the range conditions (being conditions on
zeros of the Fourier transform) survive the convolution. Thus,
functions $g_k(x,t)$ satisfy the range conditions.

According to Theorem \ref{T:Main_ball}, each of the functions
$g_k$ can be represented as $g_k=R_S f_k$ with $f_k\in
C^\infty_0(B_{1-\varepsilon})$. Then the Proposition
\ref{P:stability} shows existence of a limit $f=\lim_{k
\rightarrow \infty}f_k$ in $H_s^0(B_{1-\varepsilon})$ such that
$g=R_S f$. This finishes the proof. \qed

\section{The case of general domains}\label{S:general_domains}
Formulation of the range conditions (2b) and (3b) of Theorem
\ref{T:Main_ball} do not use explicitly that the set $S$ of
centers is a sphere, and hence that the supports of functions
under consideration are contained in the ball $B$. Only the range
conditions (4b) (the ones involving Bessel functions) explicitly
use such rotational invariance. The reader must have noticed that
in fact we proved the necessity of conditions (2b) and (3b) of
Theorem \ref{T:Main_ball}, as well as the range conditions of
Theorem \ref{T:estimates}, for arbitrary domain $D$. The only
caveat is that, as we have discussed already, the moment
conditions (2a) and (3a) should be formulated in terms of Lemma
\ref{L:moment_general} only, rather than in terms of Lemma
\ref{L:moment_ball}. In other words, we have proven some necessary
range conditions for the spherical mean transforms with centers on
the boundary of an arbitrary smooth domain. In order to formulate
them, let us introduce a number $T$ such that every sphere
centered in $\overline{D}$ and of radius at least $T$ does not
intersect the interior of $D$. We now consider the cylinder
$C=D\times[0,T]$ and formulate the conditions
\begin{equation}\label{E:Darboux_term_T}
G(x,T)=G_t(x,T)=0.
\end{equation}
The following theorem was also proven while we were proving
Theorem \ref{T:Main_ball}:
\begin{theorem}\label{T:Main_general}
Let $D\subset \RR^n$ be a bounded domain with the smooth boundary
$\Gamma$. Consider the transform
$$
g(x,t)=R_\Gamma f(x,t)=\omega^{-1}\int\limits_{|y|=1}
f(x+ty)dS(y), x\in \Gamma, t\geq 0.
$$

Then, if $f\in C^\infty_0(D)$ and $g(x,t)=R_\Gamma f$, the
following three range conditions hold.
    \begin{enumerate}
    \item
    \begin{enumerate}
    \item The moment conditions of Lemma \ref{L:moment_general} are satisfied.
    \item The solution $G(x,t)$ of the interior problem
    (\ref{E:Darboux}), (\ref{E:Darboux_bound}), (\ref{E:Darboux_term_T}) in $C$ (which always exists for $t>0$)
    satisfies the condition
    $$
    \lim\limits_{t \to 0}\int\limits_B \frac{\partial G}{\partial t}(x,t)\phi(x)dx=0
    $$
    for any eigenfunction $\phi(x)$ of the Dirichlet Laplacian in
    $D$.
    \end{enumerate}
    \item
    \begin{enumerate}
    \item The moment conditions of Lemma \ref{L:moment_general} are satisfied.
    \item Let $-\lambda^2$ be an eigenvalue of the Dirichlet Laplacian in $D$
    and $u_\lambda$ be the corresponding eigenfunction solution (\ref{E:Darboux_eigenf}).
    Then the following orthogonality condition is satisfied:
    \begin{equation}
\int\limits_{\Gamma\times [0,T]} g(x,t)\partial_\nu u_\lambda
(x,t)t^{n-1}dxdt=0.
    \end{equation}
    Here $\partial_\nu$ is the exterior normal derivative at the
    lateral boundary of $C$.
    \end{enumerate}
\item The moments
$$M_k(x)=\int\limits_0^{\infty}g(x,t)t^{2k+n-1}dt$$
extend from $x \in \Gamma$ to $x \in \mathbb R^n$ as polynomials
$Q_k(x)$ satisfying the recurrency condition
(\ref{E:Moment_relation}) and the growth estimates
$$
|Q_k(x)|\le M^k, x \in D,
$$ for some constant $M>0$.
    \end{enumerate}
    Moreover, range conditions (1) and (2) on a function $g$ are equivalent.
    \end{theorem}
One can ask whether these conditions are sufficient (together with
appropriate smoothness and support conditions on $g$) for $g$
being in the range of the transform $R_\Gamma$. The authors plan
to address this topic in another publication.

One should also note here that in the case of a general domain
$D$, the conditions (\ref{E:Darboux_term_T}) do not take into
account geometry of the domain and can be specified further.
Namely, let us define the following function in $D$:
\begin{equation}\label{E:domain_shape}
    \rho (x)=\mathop{\mbox{sup}}\limits_{y\in\partial D}|x-y|, \,
    x\in \overline{D}.
\end{equation}
Then clearly (\ref{E:Darboux_term_T}) can be replaced by
\begin{equation}\label{E:Darboux_term_domain}
    G(x,t)=0 \mbox{ for } t\geq \rho (x).
\end{equation}

\section{Proofs of some lemmas}\label{S:lemmas}

\subsection{Proof of Lemma \ref{L:PW-Bessel}}
As it has been mentioned in text, this is a known result, so we
provide a quick sketch of the proof here for reader's convenience.

Let us prove the necessity of the conditions first. Evenness is
immediate. Let us establish Paley-Wiener estimates. Consider the
natural extension of the function $g(t)$ to a radial function
$H(y)=g(|y|)$ on $\RR^n$. Due to smoothness, evenness, and
compactness of support of $g$, we see that $H$ is a smooth
function on $\RR^n$ with the support in the ball of radius $a$.
The standard Paley-Wiener theorem now claims that the
$n$-dimensional Fourier transform $\hat H (\xi)$ of $H$ is an
entire function on $\CC^n$ with Paley-Wiener estimates analogous
to (\ref{E:PW-estimates}). Since it is known that $\Phi(\lambda)$
is just the restriction of $\hat H (\xi)$ to the set $\xi=\lambda
\theta$, where $\theta$ is  a unit vector in $\RR^n$, we get the
required estimate (\ref{E:PW-estimates}).

Now we prove the sufficiency. Consider function $F(x)=\Phi(|x|)$
on $\RR^n$. Due to conditions on $\Phi$, function $F(x)$ is smooth
everywhere outside the origin and decays with all its derivatives
faster than any power of $|x|$. So, if we can establish smoothness
at the origin, this will mean that $F$ belongs to the Schwartz
class $\mathcal{S}$. Smoothness at the origin, however,
immediately follows from the radial nature of $F$ and evenness of
$\Phi$.

Thus, there exists a radial function $H$ on $\RR^n$ of the
Schwartz class, such that its Fourier transform is equal to $F$.
If now we define a function $g(t)$ such that $H(x)=g(|x|)$, then
$g$ is the function we need. We, however, need to establish that
$g$ has support in $[-a,a]$. This is equivalent to $H$ having its
support in the ball $|x|\leq a$. Consider the standard Radon
transform of $H$:
$$
u(s,\theta)=\mathcal{R} H
(s,\theta):=\int\limits_{x\cdot\theta=s}H(x)dx,|\theta|=1.
$$
Due to the projection-slice theorem
\cite{Leon_Radon,GGG2,Helg_Radon,Natt4,Natt2001}, the one
dimensional Fourier transform $\hat u (\lambda, \theta)$ of the
Radon data $u(s,\theta)$ with respect to the linear variable $s$
coincides (up to a non-zero constant factor) with the Fourier
transform of $H$ evaluated at the point $\lambda \theta$, i.e.
with $F(\lambda \theta)$. Due to the Paley-Wiener estimates we
have for $F(\lambda \theta)$ and standard $1D$ Paley-Wiener
theorem, we conclude that $u(s,\theta)$ vanishes for any $\theta$
and any $|s|>a$. Now, since $H$ is of the Schwartz class and its
Radon transform vanishes for any $|s|>a$, the "hole theorem"
\cite{Helg_Radon,Natt4,Natt2001} implies that $H(x)=0$ for
$|x|>a$. This concludes the proof of the lemma. \qed

\subsection{A Weyl transform proof of Lemma \ref{L:equiv}}
We provide here a modification of a part of the proof of part 2 of
Lemma \ref{L:equiv}. It is based on the same Weyl transformation
from Darboux equation to the wave equation, with the inhomogeneity
moved from the boundary conditions to the right hand side. So,
after the transform $G \mapsto U=\mathcal{W}G$, where $\mathcal{W}$ denotes the Weyl
transform with respect to $t$, we, as before, arrive to the
proving of evenness with respect to time of the solution of the
following problem:
\begin{equation}\label{E:altern_wave}
    \begin{cases}
    U_{tt}-\Delta_x U=-\partial^2_t \mathcal{W}(E) \mbox{ in } D\times
    [0,T]\\
    U(x,t)=0 \mbox{ for }x\in\partial D\\
    U(x,T)=U_t(x,T)=0 \mbox{ for }x\in D.
    \end{cases}
\end{equation}
Here $E(x,t)$ is a smooth function, even with respect to $t$, and
having zero of infinite order at $t=0$.

We need to establish that $U_t(x,0)=0$.

Now the argument starts to differ somewhat from what we had
before.

Let us take any Dirichlet eigenfunction $\phi=\phi_k$ in $D.$ It
will be convenient to rewrite formula (\ref{E:Weyl}) using
integration by parts and taking into account vanishing of the
integrand in a neighborhood of $\infty:$
$$
\langle U(\cdot,t), \phi\rangle=\langle
\left(\mathcal{W}G\right)(\cdot,t),\phi\rangle= \mbox{const}
\int\limits_t^\infty \langle\partial_s G(\cdot,s),\phi\rangle
(s^2-t^2)^{(n-1)/2} ds,
$$
where $\langle\cdot,\cdot\rangle$ denotes the inner product in
$L^2(D).$

Differentiating with respect to $t$ yields
\begin{equation} \label{E:U_t}
\langle U_t(\cdot,t),\phi\rangle= -\mbox{const}(n-1) t
\int\limits_t^{\infty} \langle\partial_s G(\cdot,s),\phi\rangle
(s^2-t^2)^{(n-3)/2} ds
\end{equation}
We have proven in part 1 that the conditions (a) and (b) imply
that
$$
\lim\limits_{t \to 0+}
\langle\partial_tG(\cdot,t),\phi\rangle=0.
$$
Therefore, the
integral in the right side hand tends to a finite limit as $t \to
0+$ and, taking into account the presence of the factor $t$ in
(\ref{E:U_t}), we get
$$
\langle\partial_t U(\cdot,0),\phi\rangle= \lim\limits_{t \to 0}
\langle U_t(\cdot,t),\phi\rangle=0.
$$
Thus, the function $\partial_t U(x,0)$ is orthogonal to arbitrary
Dirichlet eigenfunction $\phi(x)$ and hence
$$\partial_t U(x,0)=0, x \in D.$$

\noindent Then, due to the wave equation, all derivatives
$\partial_t^pU(x,0)$ of odd orders $p$ vanish. Hence, iterates of
the differential operator $\partial/\partial
(t^2)=(2t)^{-1}\partial_t$ preserve smoothness of $U(x,t)$ at
$t=0$, as a function with values in $H^s(D)$. Formula
(\ref{E:Weyl}) implies then that $G=\mathcal{W}^{-1}U$ is continuous and
differentiable at $t=0$. Moreover, $\partial_tG(x,0)=0$, since
$\partial_tG(\cdot,0)$ is orthogonal to all Dirichlet
eigenfunctions $\phi$. This finishes the proof of the sufficiency
part 2 of Lemma \ref{L:equiv}.

\subsection{An alternative proof of estimates on coefficients $b_l$}
We provide an alternative growth estimate derivation for the
coefficients $b_l$ (see equation (\ref{E:lambda_boundary}) and
considerations after it). We would like to establish Paley-Wiener
estimates for the functions $b_{l}(\lambda)$, $\lambda \in \CC$.
These estimates can be derived from the identity (\ref{E:lambda})
by averaging over all $0 \le r \le 1$ the point-wise estimates
that follow from (\ref{E:lambda}). Namely, let us take the
absolute value in both sides in (\ref{E:lambda}) raised to a fixed
power $\gamma > 1$ (to be specified later) and integrate with
respect to $r$ from $0$ to $1$. The triangle inequality leads to
\begin{equation}
|b_{l}(\lambda)|^{\gamma} \int\limits_0^1 |j_{n/2-1+m}(\lambda
r)|^{\gamma}|\lambda r|^{\gamma m}dr \le \int\limits_0^1
\int\limits_{\theta \in S} |\widehat
G^m(r,\lambda)|^{\gamma}|Y^{m}_{l}(\theta)|^{\gamma} dS(\theta)
dr.\label{E:growth}
\end{equation}
Since $G^m(x,t)=0$ for $t>2$, the Paley-Wiener estimate
$$
|\widehat G^m(r\theta,\lambda)| \le C_N (1+|\lambda|)^{-N} e^{2|Im
\lambda|}
$$
 holds uniformly with respect to $r \in [0,1], \theta \in S$. Then
(\ref{E:growth}) gives:
\begin{equation}\label{E:growth1}
|b_{l}(\lambda)|^{\gamma} D(\lambda) \le const_N
(1+|\lambda|)^{-\gamma N} e^{2\gamma |Im \lambda|},
\end{equation}
where
$$
D(\lambda)= \int\limits_0^1|j_{n/2-1+m}(\lambda
r)|^{\gamma}|\lambda r|^{\gamma m}dr.
$$
The integral $D(\lambda)$ satisfies the estimate
\begin{equation}\label{E:growth2}
|\lambda|D(\lambda) \ge D_0 >0.
\end{equation}
Indeed, change of variable $r=|\lambda|^{-1}z$ in the integral
$D(\lambda)$ yields
$$
D(\lambda)= |\lambda|^{-1} \int\limits_0^{|\lambda|}|z^m
j_{n/2-1+m}(\frac{\lambda}{|\lambda|}z)|^{\gamma}dz.
$$
The following estimate of Bessel functions at $\infty$ is well
known:
$$
|z^m j_{n/2-1+m}(z)|^{\gamma} \le \frac{C}{|z|^{(n-1)/2}}.
$$
Thus, if we take $\gamma>2n/(n-1)$, the integral converges:
$$
\int\limits_0^{\infty}|z^m j_{n/2-1+m}(z)|^{\gamma}dz < \infty,
$$
and therefore $|\lambda||D(\lambda)|$ tends, as $|\lambda| \to
\infty$, to a finite positive constant. This gives us the required
estimate (\ref{E:growth2}) from below.

Now substituting (\ref{E:growth2})in (\ref{E:growth1}) and raising
both sides of the inequality to the reciprocal power $1/\gamma$
yields:
$$
|b_{l}(\lambda)| \le const_N (1+|\lambda|)^{-N}e^{2|Im \lambda|}.
$$
Finally, combining this estimate for $b_{l}(\lambda)$ with
Paley-Wiener estimate for $j_{n/2-1+m}(\lambda r)$,  we obtain
from (\ref{E:FB}) the needed Paley-Wiener estimate:
\begin{equation}\label{E:PW-G}
|\widehat G^m(r\theta,\lambda)| \le const_N (1+|\lambda|)^{-N}
e^{(r+2)|Im \lambda|}.
\end{equation}

\section{Final remarks}\label{S:remarks}

\begin{itemize}

\item The results of Theorems \ref{T:Main_ball} -- \ref{T:sobolev}
easily rescale by change of variables to the ball $B$ of arbitrary
radius, which we leave as a simple exercise to the reader.

\item It is necessary to note that the range condition (4) of
Theorem \ref{T:Main_ball} is an extension of the two-dimensional
one in \cite{AmbKuch_range}. However, one notices a weaker
formulation of the moment conditions than in \cite{AmbKuch_range}.

The condition (3) is a reformulation of (4), but, unlike (4), it
is suitable for arbitrary domains.

Condition (2) is not directly the one of \cite{FinRak}, but they
are definitely of the same spirit. The authors of \cite{FinRak}
work with the wave equation, while we do with Darboux. In fact,
except the dimension $3$, we work with different (albeit related)
transforms.

\item The condition 1(b) in Lemma \ref{L:equiv} and thus 2(b) in
Theorems \ref{T:Main_ball} and \ref{T:sobolev}, as well as
condition 2 in Theorem \ref{T:Main_odd} can be replaced by much
weaker ones on the behavior of $\partial_t G(x,t)$ or $G(x,t)$ at
$t=0$. One only needs to ensure that the function $h(t)$
constructed in the proof has singularities at $t=0$ that are
milder than those of Bessel functions of the second kind, and thus
has no singularities at all. For instance, one can request that
$$
\int\limits_D\partial_tG(x,t)\phi(x)dx=o(t^{1-n}), t \to 0+.
$$
For the same reasons, the above conditions for $\partial_tG(x,t)$
can be replaced by analogous conditions for $G(x,t)$, with
$t^{1-n}$ replaced by $t^{2-n}$ for $n>2$ and $\log t$ for $n=2$.

\item It is proven in Theorem \ref{T:Main_odd} that, similarly to
the results of \cite{FinRak}, moment conditions are not needed in
odd dimensions. We suspect that one cannot remove the moment
conditions in even dimensions, albeit we do not have a convincing
argument at this point. One notices, however, that Huygens'
principle played significant role in the proof of Theorem
\ref{T:Main_odd}.

One can also notice that, as it is shown in \cite{AmbKuch_range},
moment conditions alone do not suffice for the complete
description of the range. However, Theorem \ref{T:estimates} shows
that a strengthened version of moment conditions does suffice.

\item It is instructive to note that we proved in Section
\ref{SS:global} the following statement: Let $G_+(x,t)$ be the
interior solution of the Darboux equation constructed using the
orthogonality conditions inside the cylinder $C=B\times \RR$ with
zero data at $t=2$. If $G_+(x,0)$ vanishes with all its
derivatives at the boundary of $B$, then it extends to the global
solution of the Darboux equation $G$, with the initial data
supported in $B$. This holds in any (not necessarily odd, like in
Section \ref{SS:global}) dimension. Thus, the job of the moment
conditions is to ensure smooth vanishing of the interior solution
$G_+$ at the boundary, and as we proved, this is automatic in odd
dimensions.

\item An interesting interpretation of the range conditions (a)
and (b) in Theorem \ref{T:Main_ball} comes from a model presented
in \cite{Agr_parab}, where a link is established between the
spherical mean Radon transform in $\RR^d$ and the planar Radon
transform of distributions supported on a paraboloid in
$\RR^{n+1}$. We plan to discuss this interpretation in detail
elsewhere (see also a general approach to such relations in
\cite{Gi,GGG2}).

\item A natural question is: why do we need to restrict the
support of the function $f$ to the interior of the surface
$\Gamma$ of the centers? Cannot the range of $R_\Gamma$ be
reasonably described for compactly supported functions with
supports reaching beyond the surface $\Gamma$? As it was explained
in \cite{AmbKuch_range}, this is not to be expected. Briefly,
microlocal arguments of the type the ones in \cite{LQ,Q1993,XWAK}
show that the range would not be closed in reasonable spaces
(e.g., in Sobolev scale), which is a natural precondition for
range descriptions of the kind described in this paper.

\item One notices that use of microlocal tools as in
\cite{AmbKuch_range} was essentially avoided in this text (except
of the proof of Theorem \ref{T:sobolev}). What replaced it, is
using properties of the solution of Darboux equation instead.
Namely, existence of the solution $G(x,t)$ in $C$, and especially
its regular behavior at $t=0$ did the job. Thus, microlocal
analysis was replaced by much simpler PDE tools (simple properties
of the wave equation and Fourier transforms).

\item Proposition \ref{P:stability} that concerns stability
estimates \cite{Pal_funk}, answers a question raised in the remark
section of \cite{AmbKuch_range}.

\item The range condition (2) of  Theorem \ref{T:Main_ball} also
provides a reconstruction procedure: one goes from $g(x,t)$ to the
solution $G(x,t)$ of the reverse time initial-boundary value
problem for Darboux equation in $C$, and then sets $f(x):=G(x,0)$.
A similar effect was mentioned in \cite{FinRak} concerning their
range conditions, where wave equation replaces Darboux.

Such a procedure is not that abstract. It essentially corresponds
to the time reversal reconstruction. Similar consideration lead
the authors of \cite{FPR} to their reconstruction formulas in odd
dimensions (which has been extended recently to all dimensions
\cite{FHR,Kunyansky}) and the authors of \cite{MXW2} to what they
called ``universal reconstruction formula'' in $3D$.

One can notice that the proof of our range theorem also involves
implicitly a reconstruction procedure based on eigenfunction
expansions, assuming the knowledge of the full spectrum and
eigenfunctions of the Dirichlet Laplacian in $B$. This procedure
(which is somewhat similar to the one of \cite{Norton}) would
involve at a certain stage division of analytic functions with the
denominator having zeros. When the data is in the range of the
spherical mean transform, the range theorem guarantees
cancellation of zeros. However, such a procedure would be unstable
to implement. A better version of an eigenfunction expansion
inversion procedure, which does not involve unstable divisions,
has been developed recently by L.~Kunyansky \cite{Kunyansky}. The
knowledge of the whole spectrum and eigenfunctions is available
only in rare cases (besides when the domain is a ball), e.g. for
crystallographic domains \cite{Berard,Berard2}.

\end{itemize}

\section*{Acknowledgments}
The work of the second author was supported in part by the NSF
Grants DMS-9971674, DMS-0002195, and DMS-0604778. The work of the
third author was supported in part by the NSF Grants DMS-0200788
and DMS-0456868. The authors thank the NSF for this support. Any
opinions, findings, and conclusions or recommendations expressed
in this paper are those of the authors and do not necessarily
reflect the views of the National Science Foundation.

The authors thank G.~Ambartsoumian, W.~Bray, D.~Finch,
F.~Gonzalez, A.~Greenleaf, S.~Helgason, D.~Khavinson,
V.~Kononenko, L.~Kunyansky, V.~Palamodov, S.~Patch, B.~Rubin,
K.~Trimeche, and B.~Vainberg for useful information and
discussions.

The first author would like to thank Texas A\&M University and
Tufts University for hospitality and support. The third author
expresses his appreciation to the Tufts University FRAC and the
Gelbart Institute of Bar Ilan University for their support.

\noindent agranovs@macs.biu.ac.il\\
kuchment@math.tamu.edu\\
Todd.Quinto@tufts.edu


\addcontentsline{toc}{section}{References}
\begin{thebibliography}{99}
\bibitem{Agmon_helm1} S.~Agmon, {\it A representation theorem for solutions
of the Helmholtz equation and resolvent estimates for the
Laplacian}, in: ``Analysis, et cetera'' (P.~Rabinowitz and
E.~Zehnder, eds.), 39--76, Academic Press, Boston, 1990.

\bibitem{Agmon_helm2}  S.~Agmon, {\it Representation theorems for solutions of
the Helmholtz equation on $\RR^n$ }, in: ``Differential Operators
and Spectral Theory'' (M.~Sh.~Birman's 70th anniversary
collection), AMS Translations Ser.~2, Vol.~{\bf 189}, 27--43,
Amer. Math. Soc., Providence, 1999.

\bibitem{Agr} M.~L.~Agranovsky, Radon transform on polynomial level sets and
related problems, Israel Math.Conf.Proc.,11,1997.

\bibitem{Agr_parab} M.~Agranovsky, On a problem of injectivity for the Radon
transform on a paraboloid, in Analysis, geometry, number theory:
the mathematics of Leon Ehrenpreis (Philadelphia, PA, 1998),
Contemp.Math. v. 251, AMS, Provodence, RI, 2000, 1--14.

\bibitem{ABK} M.~Agranovsky, C.~A.~Berenstein, and P.~Kuchment,
Approximation by spherical waves in $L^p$-spaces, J. Geom. Anal.,
{\bf 6}(1996), no. 3, 365--383.

\bibitem{AgrNara} M.~L.~Agranovsky, E.~K.~Narayanan, Injectivity of the spherical
mean operator on the conical manifolds of spheres, Siberian
Math.J., {\bf 45} (2004), no. 4,597--605.

\bibitem{AQ} M.~L.~Agranovsky, E.~T.~Quinto, Injectivity sets for the
Radon transform over circles and complete systems of radial
functions, J. Funct. Anal., {\bf 139} (1996), 383--413.

\bibitem{AQ2} M.~L.~Agranovsky, E.~T.~Quinto, Geometry of stationary sets for the
wave equation in $\mathbb R^n$:the case of finitely suported
initial data, Duke Math.J., {\bf 107} (2001), no. 1, 57--84.

\bibitem{AQ3} M.~L.~Agranovsky, E.~T.~Quinto, Stationary sets for the wave equation
in crystallographic domains, Trans. AMS, {\bf 355} (2003), no. 6,
2439--2451.

\bibitem{AQ4} M.~L.~Agranovsky, E.~T.~Quinto, Remarks on stationary sets for the
wave equation, Integral Geometry and Tomography, Contemp.Math., v.
405, 2006, 1--11.

\bibitem{AgrVolchZal} M.~L.~Agranovsky, V.~V.~Volchkov, L.~Zalcman, Conical
uniqueness sets for the spherical Radon transform, Bull.London
Math.Sos., {\bf 31} (1999), no. 4, 363--372.

\bibitem{Akhiezer1} N.~I.~Akhiezer, To the theory of coupled and integral equations,
Zap. Kharkov Mat. Obsch. {\bf 25} (1957), 5--31. (in Russian)

\bibitem{Akhiezer_book} N.~I.~Akhiezer, \textit{Lectures on Integral
Transforms}, AMS, Providence, RI 1988.

\bibitem{AmbKuch} G.~Ambartsoumian and P.~Kuchment, On the injectivity of
the circular Radon transform arising in thermoacoustic tomography,
Inverse Problems {\bf 21} (2005), 473--485.

\bibitem{AmbKuch_range} G.~Ambartsoumian and P.~Kuchment, A range description
for the planar circular Radon transform, SIAM J. Math. Anal. {\bf
38} 2006, no. 2, 681--692.

\bibitem{And}L.-E.~Andersson, On the determination of a function
from spherical averages, SIAM J. Math. Anal. {\bf 19} (1988), no.
1, 214--232.

\bibitem{Asgeirsson} L.~Asgeirsson, \"{U}ber eine Mittelwerteigenschaft von
L\"{o}sungen homogener linearer partieller Differentialgleichungen
zweiter Ordnung mit konstanten Koeffizienten, Ann. Math., {\bf
113} (1937), 321--346.

\bibitem{Berard} P.~B\'{e}rard, Spectres et groupes cristallographiques,
C. R. Acad. Sci. Paris Sér. A-B {\bf 288} (1979), no. 23,
A1059--A1060.

\bibitem{Berard2} P.~B\'{e}rard, G.~Besson, Spectres et groupes cristallographiques.
II. Domaines sph\'{e}riques, Ann. Inst. Fourier (Grenoble) {\bf
30} (1980), no. 3, 237--248.

\bibitem{etti} E. Bouzaglo-Burov, Inversion of spherical Radon transform, methods
and numerical experiments, MS Thesis, Bar-Ilan University,
2005,1--30.

\bibitem{Siegel} M.~Chamberland, D.~Siegel, Polynomial solutions to Dirichlet problems,
Proc. Amer. Math. Soc. {\bf 129} (2001), no. 1, 211--217.

\bibitem{CQ1} A. Cormack and E.T. Quinto, A problem in radiotherapy:
questions of non-negativity, Internat. J. Imaging Systems and
Technology, 1(1989), 120--124.

 \bibitem{CQ2} A. Cormack and E.T. Quinto, The mathematics and
physics of radiation dose planning, Contemporary Math. {\bf
113}(1990), 41-55.

\bibitem{CH} R.~Courant and D.~Hilbert, \textit{Methods of
Mathematical Physics, Volume II Partial Differential Equations},
Interscience, New York, 1962.

\bibitem{Delsarte} J.~Delsarte, Une extension nouvelle de la
th\'{e}orie des fonctions presque-p\"{e}riodiques de Bohr, Acta
Math. {\bf 69}(1938), 259-317.

\bibitem{Den}A.~Denisjuk, Integral geometry on the family of
semi-spheres. Fract. Calc. Appl. Anal.  {\bf 2}(1999), no. 1,
31--46.

\bibitem{Leon_Radon}L.~Ehrenpreis, \textit{The Universality of the Radon
Transform}, Oxford Univ. Press 2003.

\bibitem{Evans} L.~C.~Evans, \textit{Partial Differential
Equations}, AMS, Providence, RI 1998.

\bibitem{Faw}J.~A.~Fawcett, Inversion of $n$-dimensional spherical
averages, SIAM J. Appl. Math.  {\bf 45}(1985), no. 2, 336--341.

\bibitem{FHR} D.~Finch, M.~Haltmeier, and Rakesh, Inversion of spherical means
and the wave equation in even dimensions, preprint 2006.

\bibitem{FinRak} D.~Finch and Rakesh, The range of the spherical mean value
operator for functions supported in a ball, Inverse Problems {\bf
22} (2006), 923-938.

\bibitem{FPR}D.~Finch, Rakesh, and S.~Patch, Determining a function from
its mean values over a family of spheres, SIAM J. Math. Anal. {\bf
35} (2004), no. 5, 1213--1240.

\bibitem{Flatto} L.~Flatto, D.~J.~Newman, H.~S.~Shapiro, The level curves of harmonic
functions, Trans. Amer. Math. Soc. {\bf 123} (1966), 425--436.

\bibitem{GGG1} I.~Gelfand, S.~Gindikin, and M.~Graev, Integral geometry
in affine and projective spaces, J. Sov. Math. 18(1980), 39-167.

\bibitem{GGG2} I.~Gelfand, S.~Gindikin, and M.~Graev, \textit{Selected Topics
in Integral Geometry}, Transl. Math. Monogr. v. 220, Amer. Math.
Soc., Providence RI, 2003.

\bibitem{GelfVil} I.~Gelfand, M.~Graev, and N.~Vilenkin,
\textit{Generalized Functions, v. 5: Integral Geometry and
Representation Theory}, Acad. Press 1965.

\bibitem{Gi} S. Gindikin, Integral geometry on real quadrics, in
\textit{ Lie groups and Lie algebras: E. B. Dynkin's Seminar},
23--31, Amer. Math.  Soc. Transl. Ser. 2, 169, Amer. Math. Soc.,
Providence, RI, 1995.

\bibitem{Griffith} J.~L.~Griffith, On Weber transforms, J. Proc. Roy.
Soc. New South Wales {\bf 89} (1955), 109--115.

\bibitem{GU} A.~Greenleaf, G.~Uhlmann, Microlocal techniques in integral
geometry, in \textit{Integral geometry and tomography (Arcata, CA,
1989)}, 121--135, Contemp. Math., v. 113, Amer. Math. Soc.,
Providence, RI, 1990.

\bibitem{Guillemin} V. Guillemin, On some results of Gelfand in integral geometry,
proc. Symp. Pure Math. {\bf 43} (1985), 149--155.

\bibitem{GS} V. Guillemin and S. Sternberg \textit{Geometric Asymptotics},
Amer. Math. Soc., Providence, RI, 1977.

\bibitem{Guill_Ster} V. Guillemin and S. Sterenberg, Some problems of integral geometry and
some related problems in microlocal analysis,  Amer. J. Math. {\bf
101} (1979), no. 4, 915--955.

\bibitem{shapiro} L.~Hansen and H.~S.~Shapiro, Functional Equations and Harmonic
Extensions, Complex Variables, {\bf 24} (1994), 121129.


\bibitem{Helg_Radon} S. Helgason, \textit{The Radon Transform}, Birkh\"{a}user, Basel
1980.

\bibitem{Helg_groups} S.~Helgason, \textit{Groups and Geometric
Analysis}, Amer. Math. Soc., Providence, R.I. 2000.

\bibitem{Hertle2} A. Hertle, The identification problem for the
constantly attenuated Radon transform, Math. Z. 197(1988), 13-19.

\bibitem{Hor} L.~H\"{o}rmander, Fourier integral operators I, Acta Math., {\bf 127} (1971),
no. 1--2, 79–183.

\bibitem{John}F.~John, \textit{Plane Waves and Spherical Means, Applied to Partial
Differential Equations}, Dover 1971.

\bibitem{Khavinson_book}D.~Khavinson, \textit{Holomorphic Partial Differential Equations
and Classical Potential Theory}, Universidad de La Laguna Press,
1996,

\bibitem{Khavinson} D.~Khavinson and H.~S.~Shapiro. Dirichlets problem when the
data is an entire function, Bull. London Math. Soc., {\bf 24}
(1992), 456468.

\bibitem{Kipriyanov} I.~A.~Kipriyanov, \textit{Singular Elliptic Boundary Value
Problems} (in Russian), Nauka, Moscow 1997.

\bibitem{Kruger}R.~A. Kruger, P.~Liu, Y.~R.~Fang, and C.~R.~Appledorn,
Photoacoustic ultrasound (PAUS)reconstruction tomography, Med.
Phys. {\bf 22} (1995), 1605-1609.

\bibitem{K1} P. Kuchment, On positivity problems for the Radon
transform and some related transforms, Contemporary Math.,
140(1993), 87-95.

\bibitem{Kuch_AMS05} P.~Kuchment, Generalized Transforms of Radon Type and
Their Applications, in \cite{OlafQuinto}, pp. 67--91.

\bibitem{KL1} P. Kuchment and S. Lvin, Paley-Wiener theorem for
exponential Radon transform, Acta Appl. Math.  {18}(1990), 251-260

 \bibitem{KL2} P. Kuchment and S. Lvin, The range of the exponential
Radon transform, Soviet Math. Dokl.  {42}(1991), no.1, 183-184.

\bibitem{KuchQuinto} P.~Kuchment and E.~T.~Quinto, Some problems of integral
geometry arising in tomography, chapter XI in \cite{Leon_Radon}.

\bibitem{Kunyansky} L.~Kunyansky, Explicit inversion formulas for
the spherical mean Radon transform, preprint 2006.

\bibitem{Levitan} B.~M.~Levitan, Expansion in Fourier series and integrals with
Bessel functions, Uspehi Matem. Nauk (N.S.) {\bf 6} (1951). no.
2(42), 102--143. (Russian)

\bibitem{Levitan_shift} B.~M.~Levitan, \textit{Generalized Translation Operators and Some
of Their Applications}, Israel Progr. Sci. Translations, Daniel
Davey \& Co., Inc.,  Jerusalem 1964. (Translation of the 1962
Russian edition)

\bibitem{LP1}V.~Ya.~Lin and A.~Pinkus, Fundamentality of ridge
functions, J. Approx. Theory, {\bf 75} (1993), 295--311.

\bibitem{LP2}V.~Ya.~Lin and A.~Pinkus, Approximation of
multivariable functions, in {\it Advances in computational
mathematics}, H. P. Dikshit and C. A. Micchelli, eds., World Sci.
Publ., 1994, 1-9.

\bibitem{Lions_dels} J.-L.~Lions, Op\'{e}rateurs de Delsarte et
probl\`{e}mes mixtes, Bull. Soc. Math. France {\bf 84} (1950),
9--95.

\bibitem{Lions}J.-L.~Lions and E.~Madgenes, \textit{Non-Homogeneous Boundary Value Problems and
Applications}, vol.1, Springer-Verlag, Berlin-Heidelberg-New York
1972.

\bibitem{LQ}A.~K.~Louis and E.~T.~ Quinto, Local tomographic methods in
Sonar, in \textit{Surveys on solution methods for inverse
problems}, pp. 147-154, Springer, Vienna, 2000.

\bibitem{Lv} S. Lvin, Data correction and restoration in emission
tomography, pp.\  149--155 in \cite{QCK}.

\bibitem{Magnus} W.~Magnus, F.~Oberhettinger, and Raj Pal Soni,
\textit{Formulas and Theorems for the Special Functions of
Mathematical Physics}, Springer-Verlag, Berlin, New York 1966.

\bibitem{Mennes} C.~Mennesier, F.~Noo, R.~Clackdoyle, G.~Bal, and
L.~Desbat, Attenuation correction in SPECT using consistency
conditions for the exponential ray transform, Phys. Med. Biol.
{\bf 44} (1999), 2483--2510.

\bibitem{Natt3} F. Natterer, Exploiting the range of Radon transform
in tomography, in: Deuflhard P. and Hairer E. (Eds.), Numerical
treatment of inverse problems in differential and integral
equations, Birkh\"{a}user Verlag, Basel 1983.

\bibitem{Natt4}F.~Natterer, \textit{The mathematics of computerized
tomography}, Wiley, New York, 1986.

\bibitem{Natt-att} F. Natterer, Inversion of the attenuated Radon
transform, Inverse Problems 17(2001), no. 1, 113--119.

\bibitem{Natt2001}F.~Natterer and F.~W\"{u}bbeling,
\textit{Mathematical Methods in Image Reconstruction}, Monographs
on Mathematical Modeling and Computation v. 5, SIAM, Philadelphia,
PA 2001.

\bibitem{Nessibi1} M.~M.~Nessibi, L.~T.~Rachdi, K.~Trimeche, Ranges and
inversion formulas for spherical mean operator and its dual, J.
Math. Anal. Appl. 196 (1995), no. 3, 861--884.

\bibitem{Nessibi2} M.~M.~Nessibi, L.~T.~Rachdi, K.~Trimeche,
Ranges and inversion formulas for spherical mean operator and its
dual, C. R. Math. Rep. Acad. Sci. Canada 17 (1995), no. 1, 17--21.

\bibitem{Nil}S.~Nilsson, Application of fast backprojection techniques for some inverse
problems of integral geometry, Linkoeping studies in science and
technology, Dissertation 499, Dept. of Mathematics, Linkoeping
university, Linkoeping, Sweden 1997.

\bibitem{NC}C.~J.~Nolan and M.~Cheney, Synthetic aperture
inversion, Inverse Problems {\bf 18}(2002), 221--235.

\bibitem{Noo_half} F. Noo, R.~Clackdoyle, and J.--M.~ Wagner, Inversion
of the $3D$ exponential X-ray transform for a half equatorial band
and other semi-circular geometries, Phys. Med. Biol. {\bf 47}
(2002), 2727--35.

\bibitem{NW} F. Noo and J.--M. Wagner, Image reconstruction in
$2D$ SPECT with $180^{\mbox{o}}$ acquisition, Inverse Problems,
17(2001), 1357--1371.

\bibitem{Norton} S.~J.~Norton, Reconstruction
of a two-dimensional reflecting medium over a circular domain:
exact solution, J. Acoust. Soc. Am. {\bf 67} (1980), 1266-1273.

\bibitem{Nort-Linzer} S.~J.~Norton and M.~Linzer, Ultrasonic reflectivity
imaging in three dimensions: exact inverse scattering solutions
for plane, cylindrical, and spherical apertures, IEEE Transactions
on Biomedical Engineering, 28(1981), 200-202.

\bibitem{OlafQuinto} G.~Olafsson and E.~T.~Quinto (Editors), \textit{The Radon
Transform, Inverse Problems, and Tomography. American Mathematical
Society Short Course January 3--4, 2005, Atlanta, Georgia}, Proc.
Symp. Appl. Math., v. 63, AMS, RI 2006.

\bibitem{Pal}V.~P.~Palamodov, Reconstruction from limited data of arc
means, J. Fourier Anal. Appl. {\bf 6} (2000), no. 1, 25--42.

\bibitem{Pal05} V.~P.~Palamodov, Characteristic problems for the spherical mean
transform, in \textit{Complex analysis and dynamical systems II},
pp. 321--330, Contemp. Math., v. 382, Amer. Math. Soc.,
Providence, RI, 2005.

\bibitem{Pal_book} V.~P.~Palamodov, \textit{Reconstructive Integral
Geometry}, Birkh\"{a}user, Basel 2004.

\bibitem{Pal_funk} V.~P.~Palamodov, Remarks on the general Funk
transform, preprint, Tel Aviv University, August 2006.

\bibitem{Patch}S.~K.~Patch, Thermoacoustic tomography -
consistency conditions and the partial scan problem, Phys. Med.
Biol. {\bf 49} (2004), 1--11.

\bibitem{Po2} I. Ponomarev, Correction of emission tomography data.
Effects of detector displacement and non-constant sensitivity,
Inverse Problems, 10(1995)  1-8.

\bibitem{Q1980} E.~T. Quinto, The dependence of the generalized
Radon transform on defining measures, {Trans. Amer. Math. Soc.}
{257}(1980), 331--346.

 \bibitem{Quinto2} E.T. Quinto, The invertibility of rotation
invariant Radon transforms, J. Math. Anal. Appl., 91(1983),
510--522.

\bibitem{Qrange}E.~T.~Quinto, Null spaces and ranges for the classical
and spherical Radon transforms, J. Math. Anal. Appl. {\bf 90}
(1982), no. 2, 408--420.

\bibitem{Q1993} E.~T.~Quinto, Singularities of the X-ray transform
and limited data tomography in ${\mathbf R}^2$ and $\mathbf R^3$,
{SIAM J.  Math. Anal.} {24}(1993), 1215--1225.

\bibitem{Q2006} E.~T.~Quinto, An introduction to X-ray tomography
and Radon transforms, in \cite{OlafQuinto}, pp. 1--23.

\bibitem{QCK} E.T. Quinto, M. Cheney, and P. Kuchment (Editors),
\textit{ Tomography, Impedance Imaging, and Integral Geometry},
Lectures in Appl.  Math., vol. 30, AMS, Providence, RI 1994.

\bibitem{Sol3} D. Solmon, Two inverse problems for the exponential
Radon transform, in \textit{Inverse Problems in Action,}(P.S.
Sabatier, editor), 46-53, Springer Verlag, Berlin 1990.

 \bibitem{Sol4} D. Solmon, The identification problem for the
exponential Radon transform, Math. Methods in the Applied
Sciences, 18(1995), 687-695.

\bibitem{Stein} E.~M.~Stein and G.~Weiss, \textit{Introduction to Fourier analysis on Euclidean spaces},
Princeton University Press, Princeton, NJ 1971.

\bibitem{trimeche} K.~Trimeche, \textit{Generalized Harmonic Analysis and Aavelet
Packets}, Gordon \& Breach, Amsterdam 2001.

\bibitem{Vilenkin} N.~Ja.~Vilenkin, \textit{Special Functions and the Theory of Group Representations},
AMS, Providence, RI 1968.

\bibitem{MXW}M.~Xu and L.-H.~V.~Wang, Time-domain reconstruction for
thermoacoustic tomography in a spherical geometry, IEEE Trans.
Med. Imag. {\bf 21} (2002), 814-822.

\bibitem{MXW2}M.~Xu and L.-H.~V.~Wang, Universal back-projection algorithm for photoacoustic
computed tomography, Phys. Rev. E {\bf 71} (2005), 016706.

\bibitem{YXW1}Y.~Xu, D.~Feng, and L.-H.~V.~Wang, Exact frequency-domain
reconstruction for thermoacoustic tomography: I. Planar geometry,
IEEE Trans. Med. Imag. {\bf 21} (2002), 823-828.

\bibitem{YXW2}Y.~Xu, M.~Xu, and L.-H.~ V.~Wang, Exact frequency-domain
reconstruction for thermoacoustic tomography: II. Cylindrical
geometry, IEEE Trans. Med. Imag. {\bf 21} (2002), 829-833.

\bibitem{XWAK} Y.~Xu, L.~Wang, G.~Ambartsoumian, and P.~Kuchment,
Reconstructions in limited view thermoacoustic tomography, Medical
Physics 31(4) April 2004, 724-733.
\end{thebibliography}
\end{document}